\documentclass[final,onefignum,onetabnum]{siuro250211}
\usepackage{amsfonts}
\usepackage{graphicx}
\usepackage[text={6.5in,9.25in},centering]{geometry}
\usepackage[utf8]{inputenc}
\usepackage{bm}
\usepackage{amsmath}
\usepackage{natbib} 
\usepackage{geometry}
\usepackage{fancyhdr}
\usepackage{lipsum}
\usepackage{amsfonts}
\usepackage{graphicx}
\usepackage{epstopdf}
\usepackage{algorithmic}

\usepackage{enumitem}
\setlist[enumerate]{leftmargin=.5in}
\setlist[itemize]{leftmargin=.5in}

\newsiamremark{remark}{Remark}
\newsiamremark{hypothesis}{Hypothesis}
\crefname{hypothesis}{Hypothesis}{Hypotheses}
\newsiamthm{claim}{Claim}
\newsiamremark{fact}{Fact}
\crefname{fact}{Fact}{Facts}

\headers{Parameter Identifiability of RNA Dynamics}{Q. Xu}

\title{Parameter Identifiability of RNA Dynamics in PDE Transport Models of Fluorescence Recovery After Photobleaching}

\author{Qinyu Xu\thanks{Duke Kunshan University 
  (\email{qx44@duke.edu}).}}

\dedication{\small\textit{Project advisor: Dr. Veronica Ciocanel\thanks{Departments of Mathematics and Biology, Duke University (\email{veronica.ciocanel@duke.edu}).}}}

\usepackage{amsopn}



\begin{document}

\maketitle

\begin{abstract}
The transport and localization of RNA molecules, crucial for cellular function and development, involve a combination of diffusion and active transport mechanisms. Here, we are motivated by understanding the dynamics of RNA in \textit{Xenopus laevis} oocytes. Fluorescence Recovery After Photobleaching (FRAP) is an experimental technique that is widely used to investigate the dynamics of molecular movement within cells by observing the recovery of fluorescence intensity in a photobleached region over time. To advance the understanding of RNA dynamics, we develop a reaction-diffusion-advection partial differential equation (PDE) model integrating both transport and diffusion mechanisms. We propose a pipeline for identifiability analysis to assess the model's ability to uniquely determine parameter values from observed FRAP data. Based on profile likelihood analysis and reparametrization, we examine the relationship between non-identifiable parameters, which improves the robustness of parameter estimation. We find out that the identifiability of the four parameters of interest is not exactly the same in different regions of the cell. Specifically, transport velocity and diffusion coefficient are identifiable in all regions of the cell, while some combinations of binding rate and unbinding rate are found to be identifiable near the nucleus.
\end{abstract}

\section{Introduction}
Studying RNA dynamics is central to understanding how gene expression is controlled in space and time during cellular processes, and it plays an especially important role in developmental biology where spatially regulated translation directs cell fate and patterning. A classical and experimentally tractable system for this is oocyte development in the frog \textit{Xenopus laevis}. According to \cite{mowry2020using}, \textit{Xenopus laevis} undergo substantial growth during maturation, reaching a diameter of roughly 1.3 mm, and develop strong polarity characterized by the spatial distribution of numerous mRNAs and proteins. For instance, Vg1 and VegT mRNAs localize to the vegetal pole and are crucial for mesoderm and endoderm induction, respectively, thereby influencing embryonic axes specification \cite{mowry1992vegetal,zhang1996xenopus}. This subcellular mRNA patterning relies on microtubule-dependent transport, with RNA-binding proteins such as Vg1RBP/Vera mediating the recognition of localization elements and guiding transport \cite{git2002kh}.
Thus, quantitatively characterizing RNA transport and localization in \textit{Xenopus laevis} oocytes provides mechanistic insight into how molecular-level dynamics drive large-scale developmental patterning, establishing a direct link between RNA movement and cell fate determination during early embryogenesis \cite{king2005putting}.

\begin{figure}[h!]
    \centering
    \includegraphics[width=0.9\textwidth]{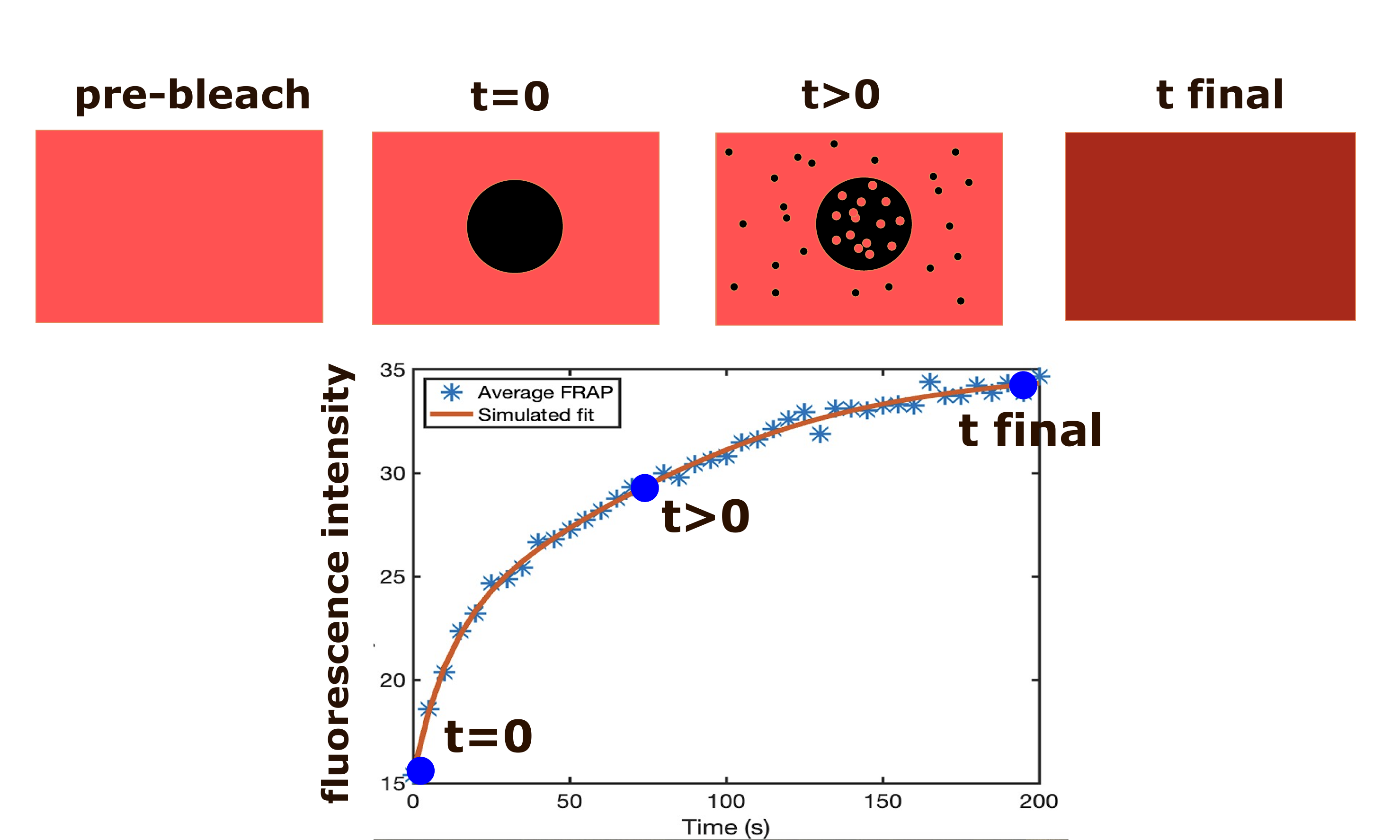}
    \caption{FRAP analysis showing: (1) Initial uniform fluorescence, (2) Sharp intensity drop after photobleaching, (3) Gradual recovery as molecules diffuse into the bleached region, and (4) Final equilibrium fluorescence. The recovery kinetics can reveal molecular mobility and binding dynamics.}
    \label{fig:frap}
\end{figure}

In living cells, RNA dynamics are commonly studied with fluorescence microscopy techniques, among which fluorescence recovery after photobleaching (FRAP) is a widely used, robust method in biological analysis. It is a standard tool in cell biology for probing protein and RNA–protein complex dynamics in their native context \cite{lippincott2001studying,reits2001fixed}. In a FRAP experiment, fluorescently tagged molecules, such as RNA–protein complexes labeled with a fluorescent protein, are visualized under a widefield microscope. A small, well-defined region of the cell is then exposed to a high-intensity laser pulse, which irreversibly photobleaches the fluorophores in that region, causing their fluorescence to drop sharply. Following bleaching, the fluorescence in the region is monitored over time as unbleached molecules from surrounding areas move in. The resulting fluorescence recovery curve reflects the mobility, transport, and exchange rates of the molecules \cite{sprague2005frap}. This process is shown in \Cref{fig:frap}.

To interpret FRAP data and relate observed recovery curves to underlying physical processes, mathematical models have been developed in order to simulate the moving behaviors of RNA inside the cells, which involve several parameters such as diffusion coefficient and transport speed. Accurate parameter estimation allows researchers to quantitatively characterize RNA behavior, validate the model against experimental data, and gain insights into the underlying biological mechanisms. Crucially, however, the ability to infer parameters depends not only on the data quality but also on identifiability — whether the model structure and available data observations allow unique recovery of parameter values. Performing identifiability analysis is essential for the accurate parameter estimation, as it reveals which parameters can be uniquely determined based on the model and observed data, while other non-identifiable parameters cannot.

This paper focuses on the identifiability analysis of parameters in a reaction–diffusion–advection model of RNA localization during \textit{Xenopus laevis} oocyte development. Building on prior studies that investigated the identifiability of the parameters in a reaction-diffusion model for RNA-binding protein PTBP3 in \textit{Xenopus laevis} oocytes \cite{ciocanel2024parameter} and the reaction-diffusion-advection model exploration for RNA dynamics \cite{ciocanel2017analysis}, we explicitly incorporate active transport into the inference framework and examine identifiability across spatially distinct oocyte regions. Our goal is to develop mathematical tools in order to fill in the gaps in the current understanding of RNA movement.

\begin{figure}[h!]
    \centering
    \includegraphics[width=0.8\textwidth]{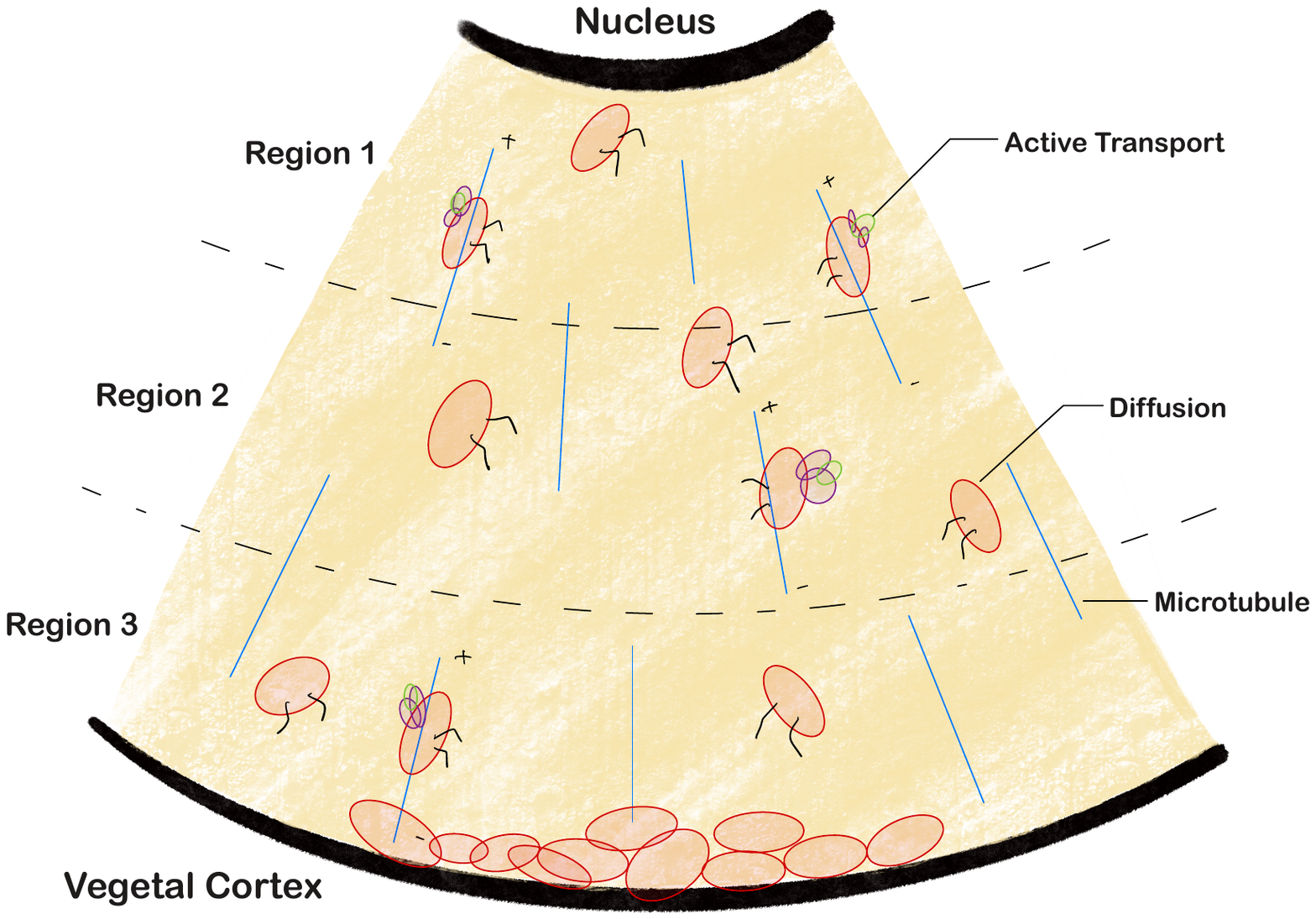}
    \caption{This cartoon shows that experiments are done in three regions of the cell based on their distance from the nucleus (top) and vegetal cortex (bottom). Red particles represent RNAs, which are transported either by active transport, the directed movement along blue microtubule tracts driven by motor proteins (purple and green blobs), or by simple diffusion. Region 1 (top) is 15 $\mu m$ from the nucleus, Region 2 (middle) is 50 $\mu m$ from the nucleus, and Region 3 (bottom) is 25 $\mu m$ from the vegetal cortex.}
    \label{fig:3regions}
\end{figure}

In the cellular environment, the nucleus is typically centrally located, while the vegetal cortex, situated opposite the animal pole, is a specialized region near the cell membrane. According to \cite{ciocanel2017analysis}, it is believed that the moving behaviors of RNA are not the same in different regions of the cell, so we divide the cell into three regions as depicted in \Cref{fig:3regions}, where Region 1 is 15 $\mu m$ from the nucleus, Region 2 is 50 $\mu m$ from the nucleus, and Region 3 is 25 $\mu m$ from the vegetal cortex \cite{ciocanel2017analysis}. For each region, a 5-$\mu m$ circular area of interest is bleached \cite{ciocanel2017analysis}, so we can obtain three average FRAP data sets from experiments. As in previous work \cite{ciocanel2017analysis}, after obtaining the baseline parameter sets that fit the real experimental FRAP data based on deterministic parameter estimation, we generate the corresponding synthetic FRAP data for the three regions. Then, the method of profile likelihood gives us the framework for determining the identifiability of parameters. We also derive the threshold for 95\% confidence interval, which serves as the criterion to decide whether a certain parameter is identifiable or not. Moreover, we investigate the parameter relationships between non-identifiable parameters, as well as apply this methodology to real and synthetic FRAP data. 

By explicitly integrating active transport into the identifiability analysis and by examining the FRAP data across biologically meaningful oocyte regions, our work extends previous efforts and provides a practical framework for quantifying which transport model parameters are identifiable when analyzing mRNA localization in development using FRAP data.

\section{Mathematical modeling of the FRAP experiment}
\subsection{Reaction-diffusion-advection model for RNA dynamics}
\begin{figure}[h!]
    \centering
    \includegraphics[width=0.6\textwidth]{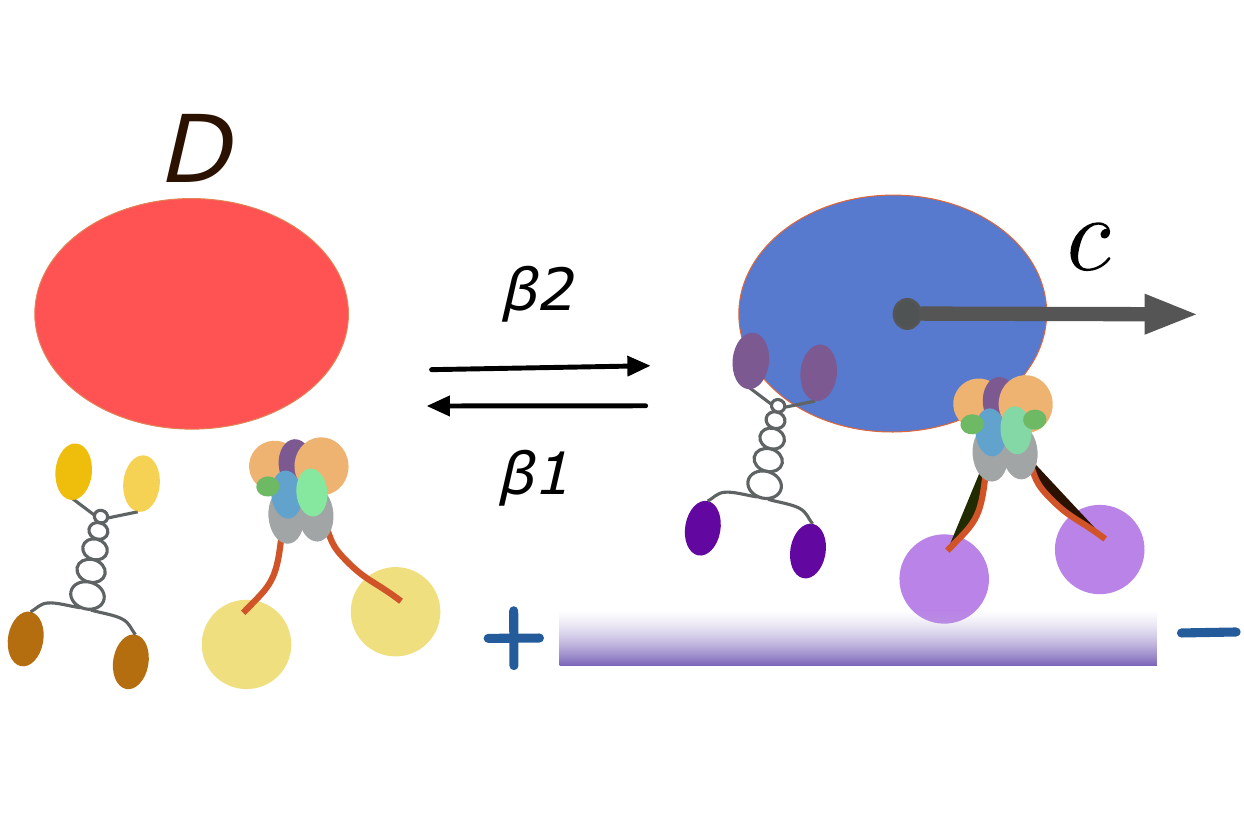}
    \caption{Cartoon of the active transport of RNA, consisting of a population of diffusing particles with diffusion coefficient $D$, a population of moving particles with velocity $c$, as well as switching rates $\beta_1$ and $\beta_2$ between the two population. In the moving state, we assume that RNA molecules are attached to motor proteins and microtubules, while in the diffusion state, we assume they are detached from microtubules.}
    \label{fig:cartoon}
\end{figure}

We model the dynamics of RNA based on two localization mechanisms. RNA can either diffuse with the diffusion coefficient $D$ in the cytoplasm, or can be transported along microtubules by motor proteins with the speed $c$. Moreover, we assume that the RNA molecules can switch between these two mechanisms with certain rates. The RNA dynamics is therefore shown in \Cref{fig:cartoon}, which is given by this reaction-diffusion-advection model \cite{ciocanel2017analysis}:

\begin{equation}
\begin{aligned}
\frac{\partial u}{\partial t} &= c \frac{\partial u}{\partial y} - \beta_1 u + \beta_2 v, \\
\frac{\partial v}{\partial t} &= D \Delta v + \beta_1 u - \beta_2 v.
\end{aligned}
\label{eq: PDEs}
\end{equation}

In this PDE model, $u$ refers to the concentration of RNA which is transported along microtubules, while $v$ refers to the concentration of RNA which is diffusing. Since the bleached area in the FRAP process is very small, we assume that within this small area, all the directions of microtubules are vertical, which is the reason why we just include $\frac{\partial u}{\partial y}$ in \Cref{eq: PDEs} (i.e., 1-dimensional transport). Furthermore, we define $\beta _1$ as the switching rate from transport RNA to diffusing RNA, which represents the unbinding rate between RNA and microtubules. Similarly, $\beta_2$ is defined as the binding rate, which is the switching rate from diffusing RNA to transport RNA. Based on this model, we are able to simulate the FRAP process using numerical PDE methods as in \cite{kassam2005fourth,cox2002exponential}. The observed FRAP recovery signal $F(t)$ is modeled as the total fluorescence within the bleached region $\Omega$, which depends on both the transport and diffusing RNA populations. Specifically, the fluorescence intensity is proportional to the sum of the concentrations $u$ and $v$, integrated over the bleach spot:
\begin{equation}
    F(t) = \int_{\Omega} \left[ u(\mathbf{x}, t) + v(\mathbf{x}, t) \right]  d\mathbf{x}
    \label{eq: PL}
\end{equation}

where $\Omega$ is the circular bleached region of diameter 5 $\mu m$, and $\mathbf{x}$ represents the spatial dimensions $(x,y)$. This expression captures the total amount of fluorescent molecules within the region of interest over time.

\subsection{Experimental FRAP data}
The experimental FRAP data used in this study, shown as the blue points in \Cref{fig:combined1}, were collected from \textit{Xenopus laevis} oocytes following the procedure described in the FRAP acquisition section of \cite{ciocanel2017analysis}. In the experiments, a small circular region ($5~\mu\mathrm{m}$ in diameter) inside the oocyte was bleached with a laser, and the fluorescence recovery in this region was recorded over 200s at regular time intervals. The recovery curves reflect how labeled RNA molecules move back into the bleached area over time.

For each region of interest, data was collected from five oocytes and their recovery curves were averaged. The fluorescence values were corrected for background changes and normalized to a nearby unbleached region to reduce experimental noise, following the adjustment procedure outlined in \cite{ciocanel2017analysis}. The resulting averaged FRAP data captures the overall recovery trend and will serve as the input for our parameter estimation procedure. In \Cref{Deterministic parameter estimation}, we fit this experimental data to \Cref{eq: PDEs} in order to extract parameter values and generate synthetic curves, which clearly shows the comparison between the model-generated synthetic curves and the measured data.


\subsection{Deterministic parameter estimation}
\label{Deterministic parameter estimation}
For the real FRAP data observed from 5 oocytes in each region, deterministic parameter estimation was performed following the numerical fitting framework in \cite{ciocanel2017analysis}. In brief, the reaction–diffusion–advection PDE model (\Cref{eq: PDEs}) was integrated numerically using an exponential time-differencing Runge–Kutta scheme with Fourier spectral discretization. The observed FRAP signal was modeled as the sum of all particle states integrated over the bleach spot, and parameters $c, D, \beta_1$ and $\beta_2$ were estimated by minimizing the $L^2$-norm between model output and experimental curves using the MATLAB routine \texttt{lsqnonlin}. The resulting best-fit parameters for each region are summarized in \Cref{tab:real FRAP}.

\begin{table}[htbp]
\footnotesize
\caption{Table of parameter estimates for the three regions \cite{ciocanel2017analysis}.}\label{tab:real FRAP}
\begin{center}
\begin{tabular}{ccccc}
\hline
\textbf{Region} & \textbf{c} & \textbf{D} & \textbf{$\beta_1$} & \textbf{$\beta_2$} \\ \hline
1 & 0.049121 & 0.258205 & 2.35E-14 & 0.006331 \\ 
2 & 0.094322 & 1.423721 & 0.003018 & 0.000762 \\ 
3 & 0.067619 & 0.830449 & 4.05E-05 & 1.37E-06 \\ \hline
\end{tabular}
\end{center}
\end{table}

\begin{figure}[htbp]
    \centering
    \begin{minipage}[t]{0.485\linewidth}
        \centering
        \includegraphics[width=\linewidth]{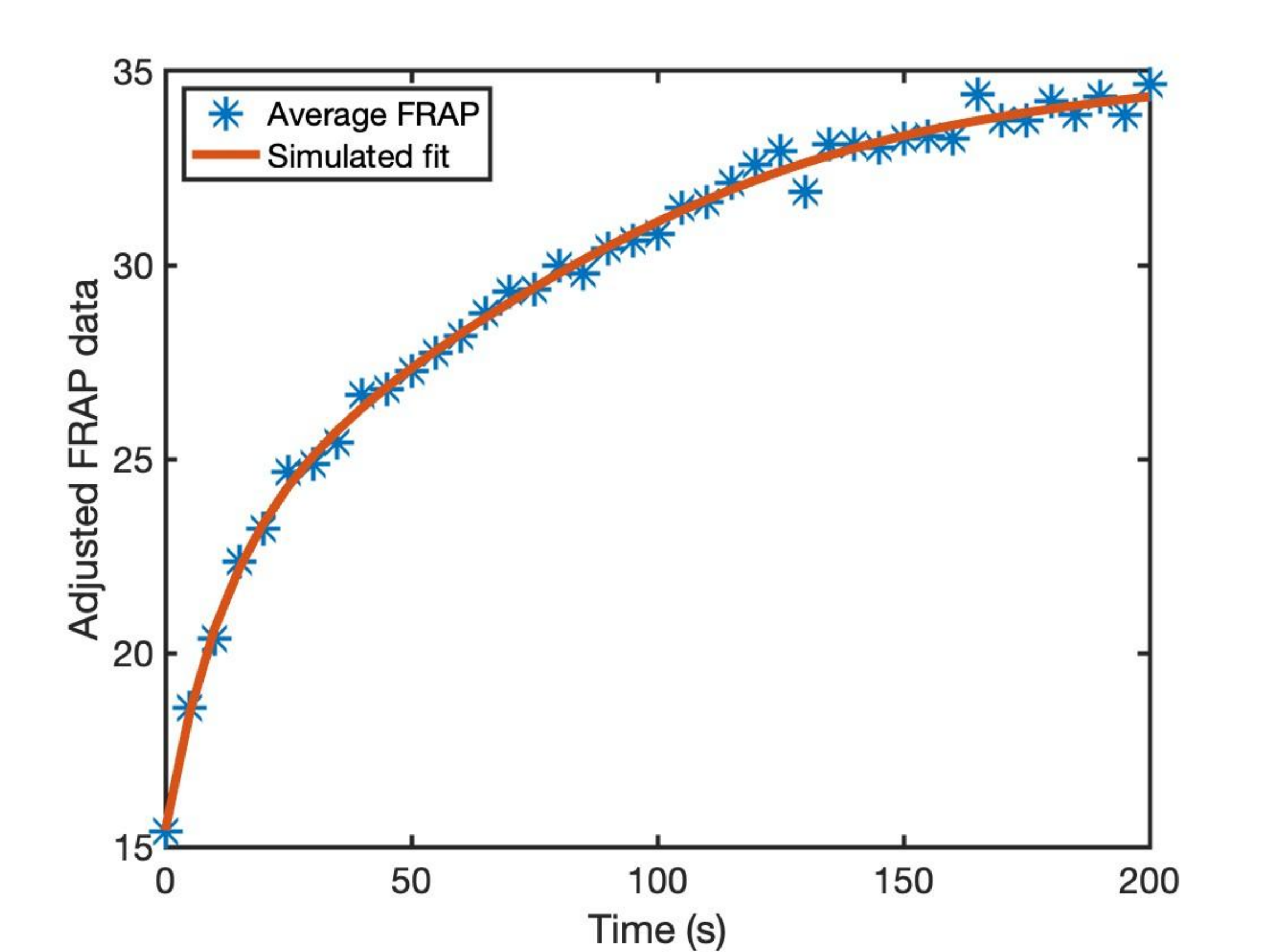}
    \end{minipage}
    \hspace{0.01\linewidth}
    \begin{minipage}[t]{0.485\linewidth}
        \centering
        \includegraphics[width=\linewidth]{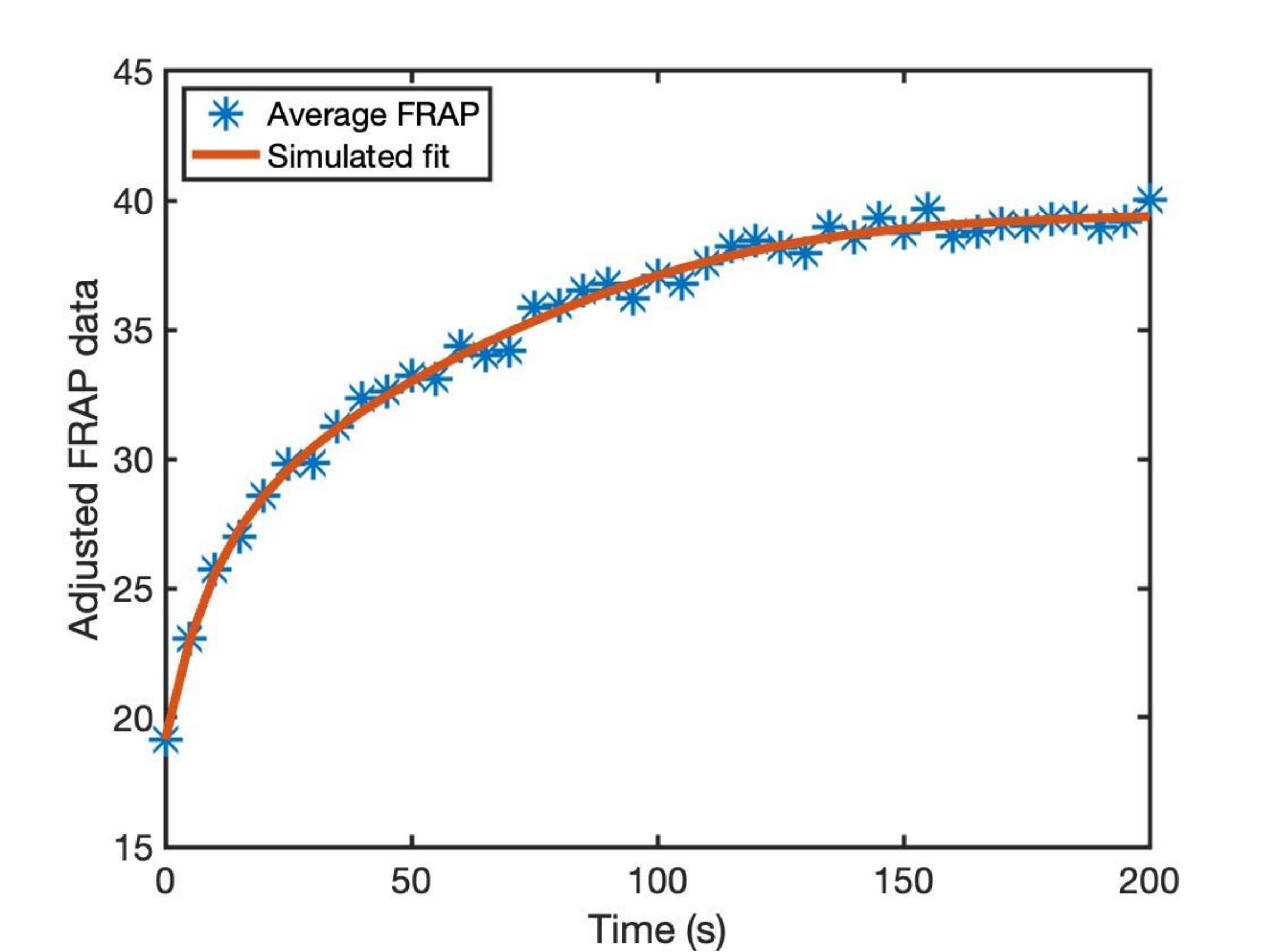}
    \end{minipage}
    
    \vspace{0.5em}
    \begin{minipage}[t]{0.485\linewidth}
        \centering
        \includegraphics[width=\linewidth]{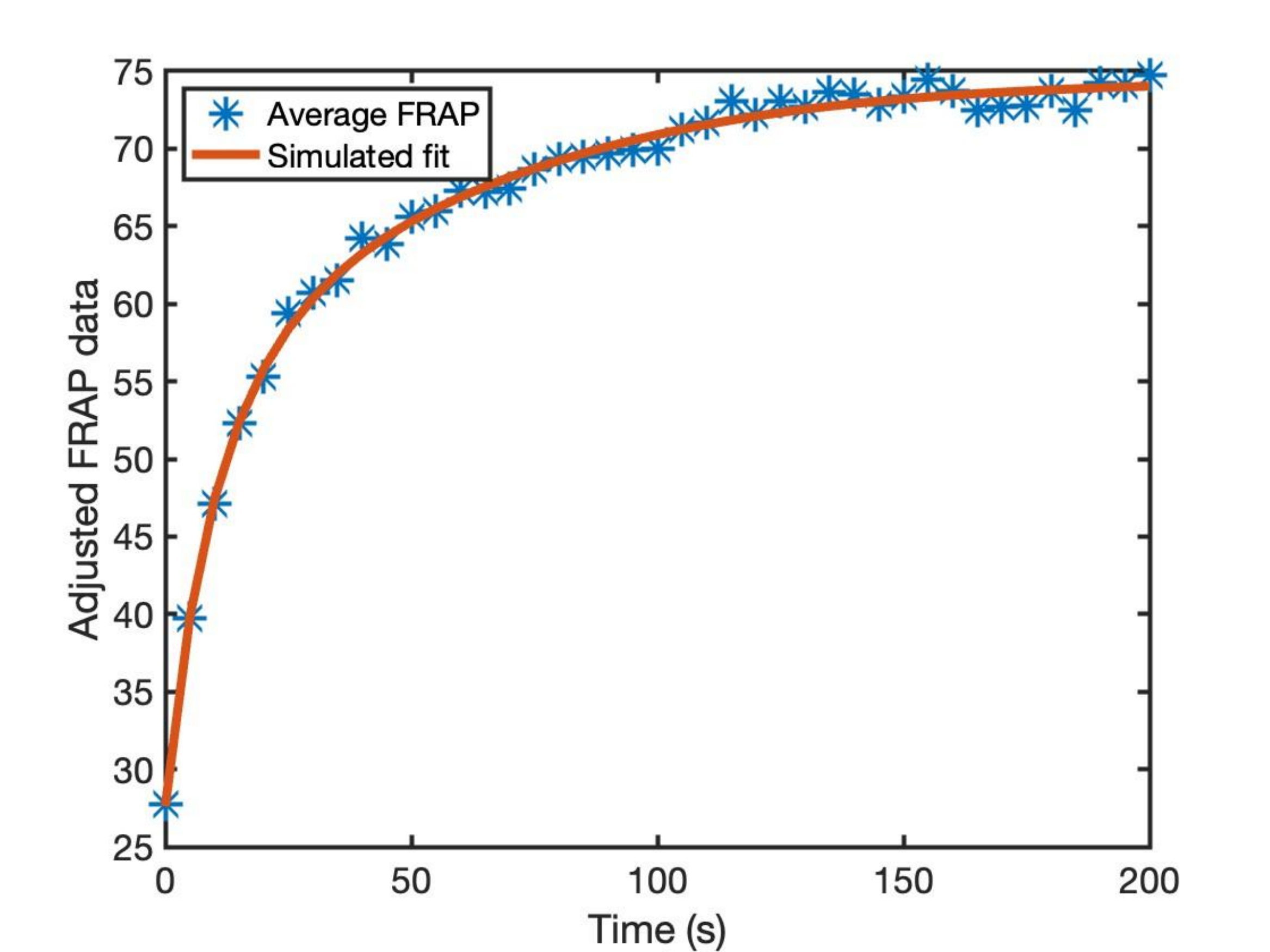}
    \end{minipage}
    
    \caption{Real FRAP data (blue points) and synthetic FRAP data (red lines) based on the baseline parameter values in \Cref{tab:baseline} (top-left: Region 1, top-right: Region 2, bottom: Region 3).}
    \label{fig:combined1}
\end{figure}

We use the parameter estimates in \Cref{tab:real FRAP} to inform some baseline parameter regimes for generating synthetic data, which are shown in \Cref{tab:baseline}. Based on these baseline values (\Cref{tab:baseline}) and synthetic FRAP data (red curves in \Cref{fig:combined1}), we are able to test our identifiability analysis framework on PDE-generated data before attempting it on real and noisy datasets \cite{ciocanel2024parameter}. It can be seen from \Cref{fig:combined1} that the red synthetic data curves fit the real experimental FRAP data points relatively well.

\begin{table}[h!]
\footnotesize
\caption{Table of baseline values for the three regions used for synthetic data generation in this work.}\label{tab:baseline}
\begin{center}
\begin{tabular}{ccccc}
\hline
\textbf{Region} & \textbf{c} & \textbf{D} & \textbf{$\beta_1$} & \textbf{$\beta_2$} \\
\hline
1 & 0.05 & 0.25 & $10^{-6}$ & $10^{-2}$ \\ 
2 & 0.1 & 1.5 & $10^{-3}$ & $10^{-4}$ \\ 
3 & 0.1 & 0.8 & $10^{-5}$ & $10^{-6}$ \\ \hline
\end{tabular}
\end{center}
\end{table}


\section{Identifiability analysis for transport models of FRAP experiments}
\subsection{Profile likelihood analysis}
We then carry out identifiability analysis based on the method of profile likelihood, starting with the synthetic datasets in \Cref{fig:combined1}. Profile likelihood analysis can provide insights into the identifiability of parameters in the model by systematically setting one parameter constant and optimizing the rest of the parameters in fitting the data. In this work, we define $\bm{\theta}$ to be the vector of all the parameters of interest: $\bm{\theta} = (c, D, \beta_1, \beta_2)$. We first choose one interest parameter $\bm{\psi}$ (for instance, $D$). By changing the values of that interest parameter, we want to find the optimized values for other nuisance parameters $\bm{\lambda}$ (for instance, $c, \beta_1, \beta_2$), so that the profile likelihood function is maximized. For instance, if we set $D$ to be the interest parameter, then by changing its value, we want to find the optimal $\bm{\lambda} = (c, \beta_1, \beta_2)$ so that the profile likelihood function is maximized. Moreover, if the measurement noise is assumed to be normally distributed, then the profile likelihood function can be expressed as:
\begin{equation}
    p(y;\bm{\psi}, \bm{\lambda}) = \left(\frac{1}{2\pi\sigma^2}\right)^{\frac{1}{2}}exp\left(-\frac{1}{2\sigma^2}||y-y_{sim}(\bm{\psi}, \bm{\lambda})||_2^2\right)
    \label{eq: PL}
\end{equation}
where $y$ refers to the noiseless synthetic data curve generated with baseline values in \Cref{tab:baseline},  $y_{sim}(\bm{\psi}, \bm{\lambda})$ refers to the noiseless synthetic data curve generated under each value of the interest parameter $\bm{\psi}$ together with the optimal values of the nuisance parameters $\bm{\lambda}$. Moreover, $||y-y_{sim}(\bm{\psi}, \bm{\lambda})||_2^2$ represents the least square error, and $\sigma$ is the standard deviation of the noise in the FRAP data. We estimate $\sigma$ by calculating the difference between each point on the synthetic data curve and the FRAP data curve generated by parameters in \Cref{tab:real FRAP}:

\begin{equation}
    \hat{\sigma} = \sqrt{\frac{\sum(y-y_{sim}(\bm{\psi}, \bm{\lambda}))^2}{N}}
    \label{eq: sigma}
\end{equation}
where $N$ is the number of data points.

\begin{figure}[h!]
    \centering
    \includegraphics[width=0.8\textwidth]{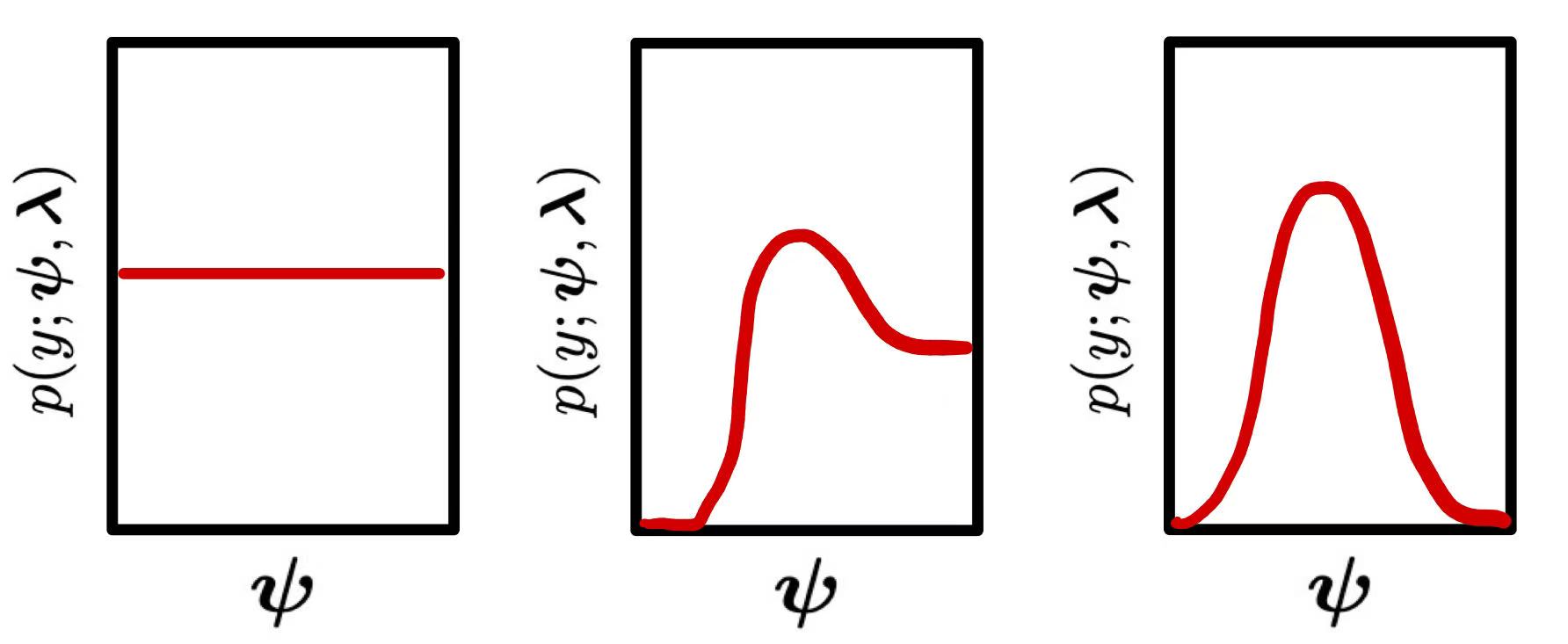}
    \caption{Interpretation of profile likelihoods. A flat likelihood (left) corresponds to structural non-identifiability, a profile that does not decrease to 0 on one or both sides of the maximum (center) indicates practical non-identifiability, and a profile with a fast decrease to 0 on both sides of the maximum (right) shows both structural and practical identifiability.}
    \label{fig:PL}
\end{figure}

For each region, we have computed the corresponding estimate values for $\sigma$:
$\hat{\sigma_1} = 0.275, \hat{\sigma_2} = 0.365, \hat{\sigma_3} = 0.614$.
\newline
To visualize the results of the profile likelihood analysis, we plot the likelihood value for each value of the interest parameter. The possible results are shown in \Cref{fig:PL} \cite{ciocanel2024parameter}. In the first plot where the curve is just a flat line, this parameter is considered as structurally non-identifiable, since the profile likelihood is not sensitive to the change of that parameter. For the second plot where the curve does not decrease to 0 on the right end, this parameter is called practically non-identifiable since its non-identifiability is mainly due to the lack of data. As for the third case, we consider this parameter to be both structurally and practically identifiable.

\subsection{Threshold for identifiability using profile likelihood analysis}
In order to establish a clear criterion to determine whether a parameter is identifiable or not using profile likelihood analysis, we calculate the threshold for the 95\% confidence interval for each profile likelihood plot. Specifically, confidence intervals are determined by applying a threshold to the likelihood ratio statistic, corresponding to a $\chi^2$ quantile. Parameters are considered identifiable when likelihood plots cross this threshold on both sides \cite{raue2009structural}. The basic formula for computing the value of this confidence threshold is given by \cite{raue2009structural,kreutz2013profile}:
\begin{equation}
    CI_{j,\alpha}(y) = \{\alpha_{p}| 2(LL(y|\bm{\hat{\theta}})-PL_{j}(p))<\Delta\alpha\}
    \label{eq: Basic CI}
\end{equation}
where $LL(y|\bm{\hat{\theta}})$ and $PL_{j}(p)$ refer to the log likelihood with the interest parameter set as value $p$:
\begin{equation}
    LL(y|\bm{\theta}) = -\frac{1}{2} \ln(2\pi\sigma^2) - \frac{\|y - y_{sim}(\bm{\psi}, \bm{\lambda})\|^2}{2\sigma^2}
    \label{eq: CI1}
\end{equation}

and $\Delta\alpha$ comes from the chi-square value ($\Delta\alpha = 3.841$ for the 95\% confidence interval). Moreover, $\bm{\hat{\theta}}$ refers to the set of parameter values that gives the maximum value of the profile likelihood. Therefore, $\bm{\hat{\theta}}$ should be the set of baseline values, and we obtain $LL(y|\bm{\hat{\theta}}) = -\frac{1}{2} \ln(2\pi\sigma^2)$ with the least square error equal to 0. Then, we plug the log likelihood into \Cref{eq: Basic CI} to get:
\begin{equation}
2(LL(y|\bm{\hat{\theta}})-PL_{j}(p)) < \Delta(\alpha)
    \label{eq: CI2.1}
\end{equation}

\begin{equation}
   \Rightarrow 2(-\frac{1}{2}ln(2\pi\sigma^2)+\frac{1}{2}ln(2\pi\sigma^2)+\frac{\|y - y_{sim}(\bm{\psi}, \bm{\lambda})\|^2}{2\sigma^2}) < \Delta(\alpha)  
    \label{eq: CI2.2}
\end{equation}

\begin{equation}
   \Rightarrow \frac{\|y - y_{sim}(\bm{\psi}, \bm{\lambda})\|^2}{\sigma^2} < \Delta(\alpha)
    \label{eq: CI2}
\end{equation}
Based on \Cref{eq: PL} for profile likelihood, the value for the threshold is given by:
\begin{equation}
    p(y;\bm{\psi}, \bm{\lambda}) >  \frac{1}{\sqrt{2\pi\sigma^2}}e^{-\frac{\Delta{\alpha}}{2}}
    \label{eq: CI3}
\end{equation}
\subsection{Outcomes of profile likelihood analysis}
\begin{figure}[htbp]
    \centering
    \begin{minipage}[t]{0.485\linewidth}
        \centering
        \includegraphics[width=\linewidth]{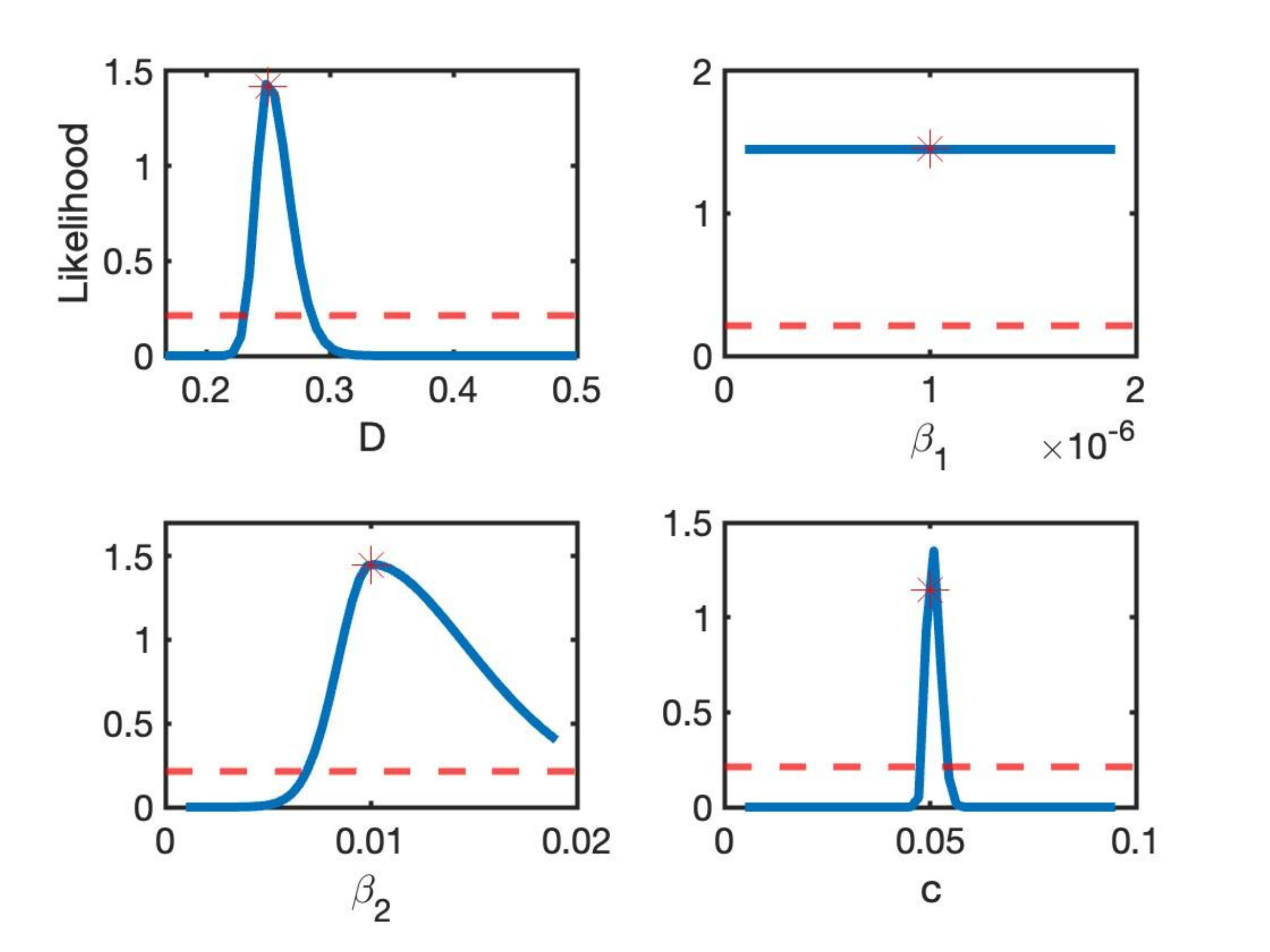}
    \end{minipage}
    \hspace{0.01\linewidth}
    \begin{minipage}[t]{0.485\linewidth}
        \centering
        \includegraphics[width=\linewidth]{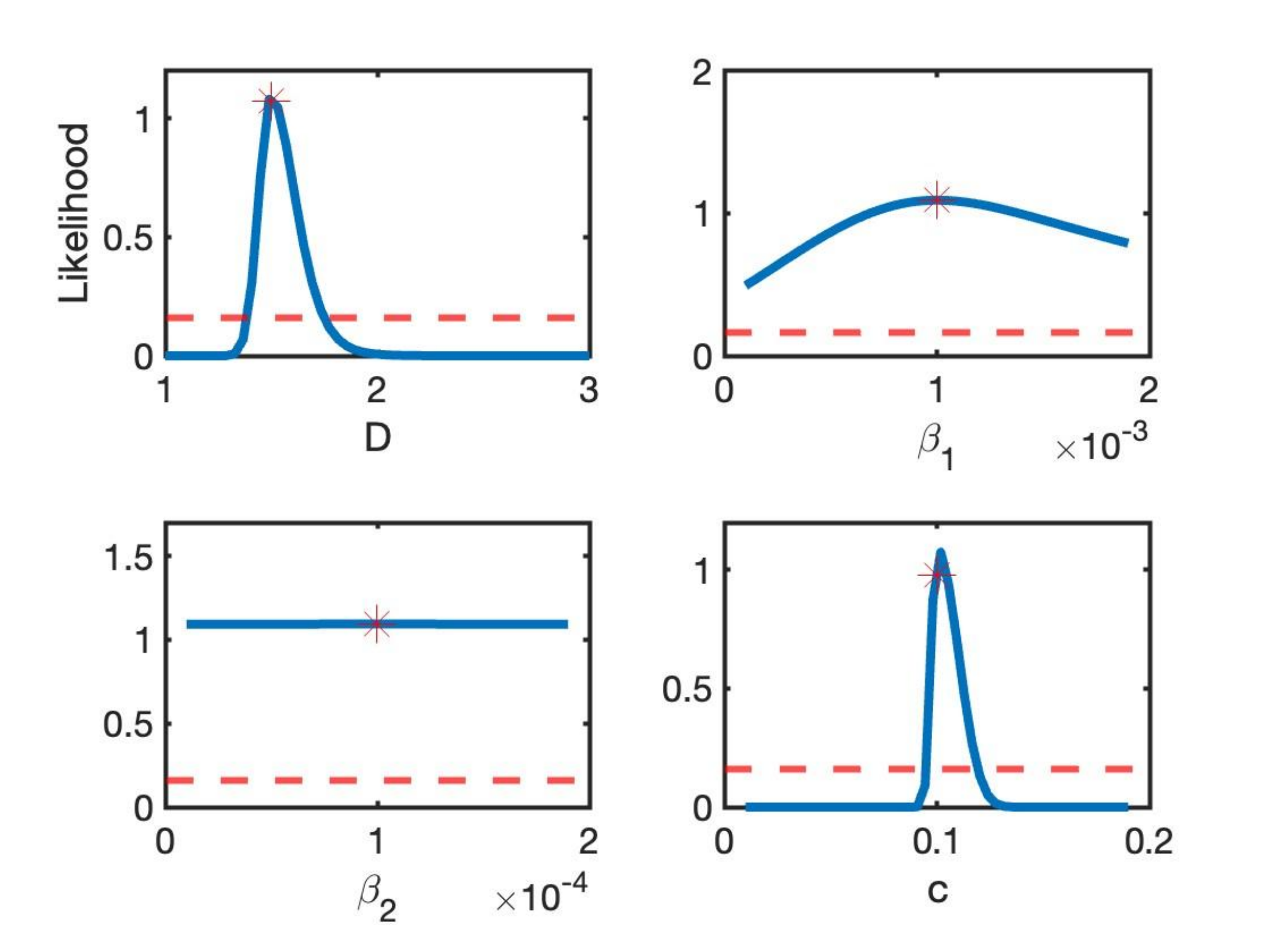}
    \end{minipage}
    
    \vspace{0.5em}
    \begin{minipage}[t]{0.5\linewidth}
        \centering
        \includegraphics[width=\linewidth]{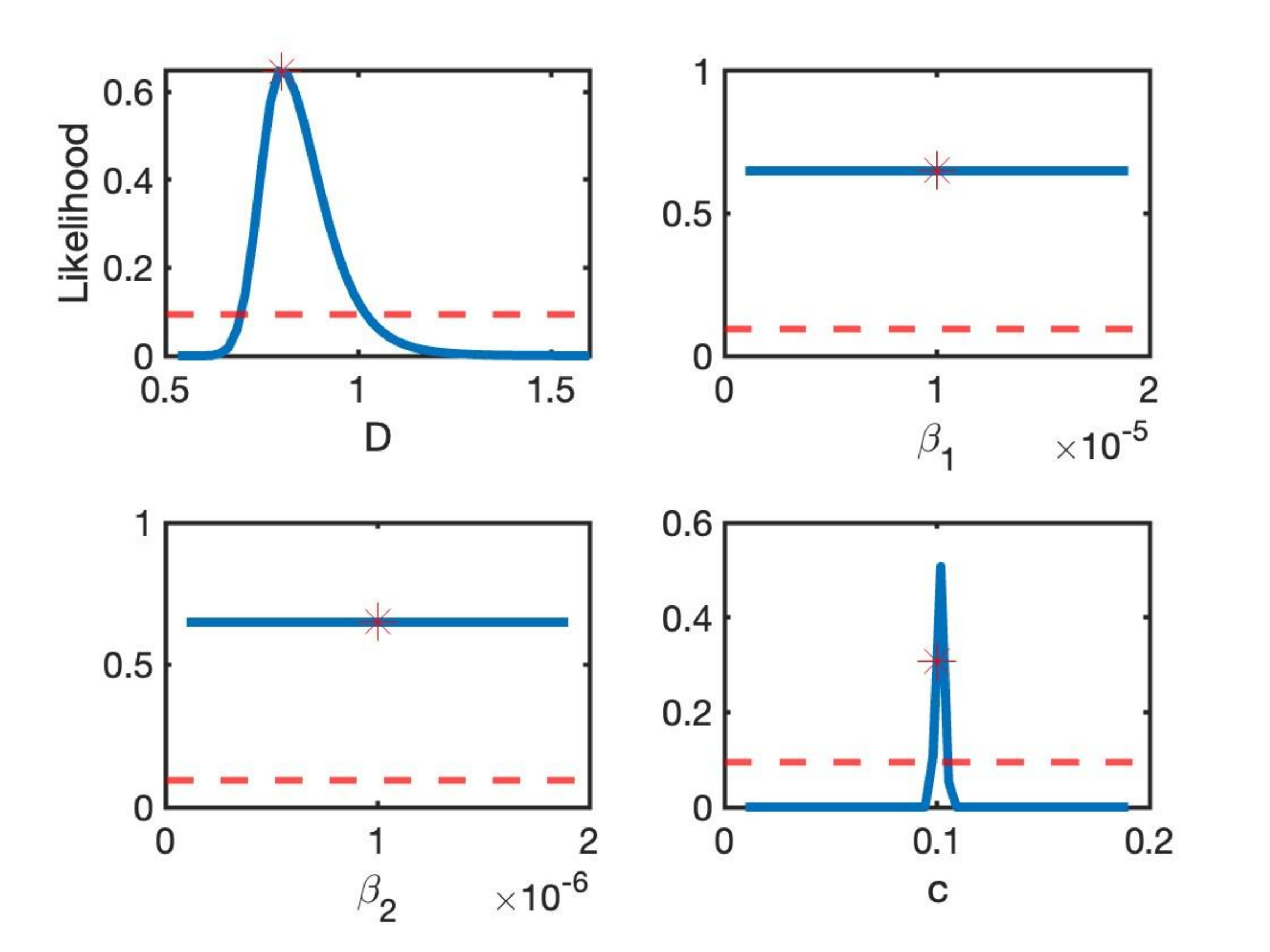}
    \end{minipage}
    
    \caption{Profile likelihood results for the three regions given noiseless FRAP data synthetically generated using \Cref{eq: PDEs} and the parameter sets in \Cref{tab:baseline}. Baseline values are indicated as red stars in each plot (top-left: Region~1, top-right: Region~2, bottom: Region~3).}
    \label{fig:PL_all}
\end{figure}

After computing the profile likelihood for all of the four parameters of interest in the three regions, and calculating the threshold for the likelihood in each plot, the results are given in \Cref{fig:PL_all}.


Overall, all three regions show similar identifiability properties for the four parameters. Notably, $c$ and $D$ are identifiable, while $\beta_1$ and $\beta_2$ are non-identifiable. Moreover, we notice that $\beta_2$ in Region 1 and $\beta_1$ in Region 2 are practically non-identifiable. As for region 3, the plots for $\beta_1$ and $\beta_2$ are just flat lines, suggesting structural non-identifiability for both kinetic rates.

\subsection{Two-dimensional profile likelihood}
In systems biology, models are becoming increasingly complex as the field strives to provide holistic descriptions of biological systems that capture both their static properties and dynamic interactions. The two-dimensional profile likelihood approach has emerged as a powerful tool by offering enhanced insights into parameter uncertainty and improving model identifiability \cite{litwin2022optimal}.

To extend the previous one-dimensional profile likelihood to a two-dimensional case, we set $c$ and $D$ to be our two interest parameters. Then, given each possible pair of $c$ and $D$ chosen from a grid around their baseline values, we want to find the optimized values for $\beta_1$ and $\beta_2$ so that the profile likelihood function is maximized.

\begin{equation}
    PL_{c, D}(p_1, p_2) = \max_{\theta \in \{\theta | c = p_1, D = p_2\}} LL(y | \bm{\hat{\theta}})
    \label{eq: 2D PL}
\end{equation}

\begin{figure}[htbp]
    \centering
    \begin{minipage}[t]{0.485\linewidth}
        \centering
        \includegraphics[width=\linewidth]{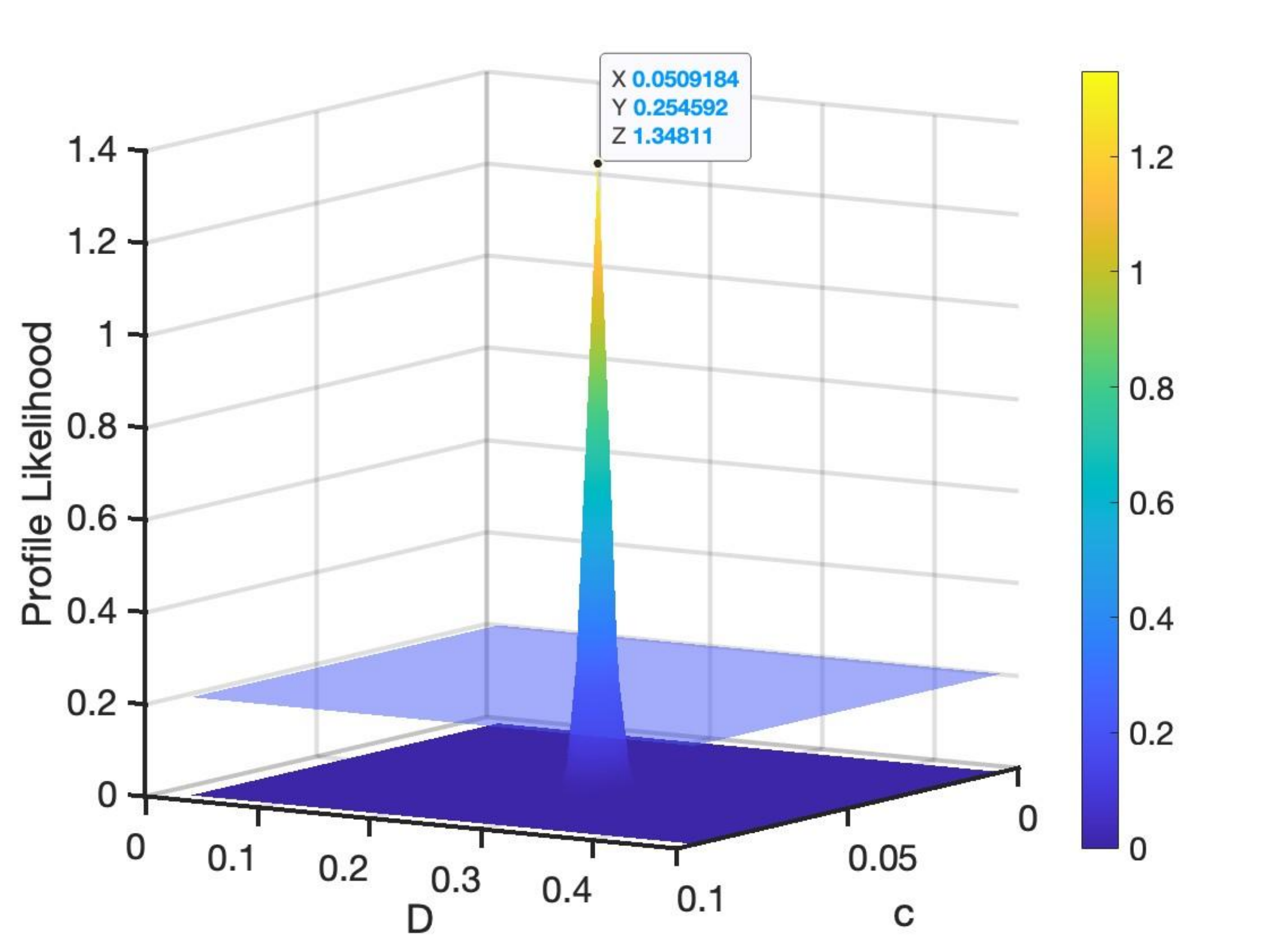}
    \end{minipage}
    \hspace{0.01\linewidth}
    \begin{minipage}[t]{0.485\linewidth}
        \centering
        \includegraphics[width=\linewidth]{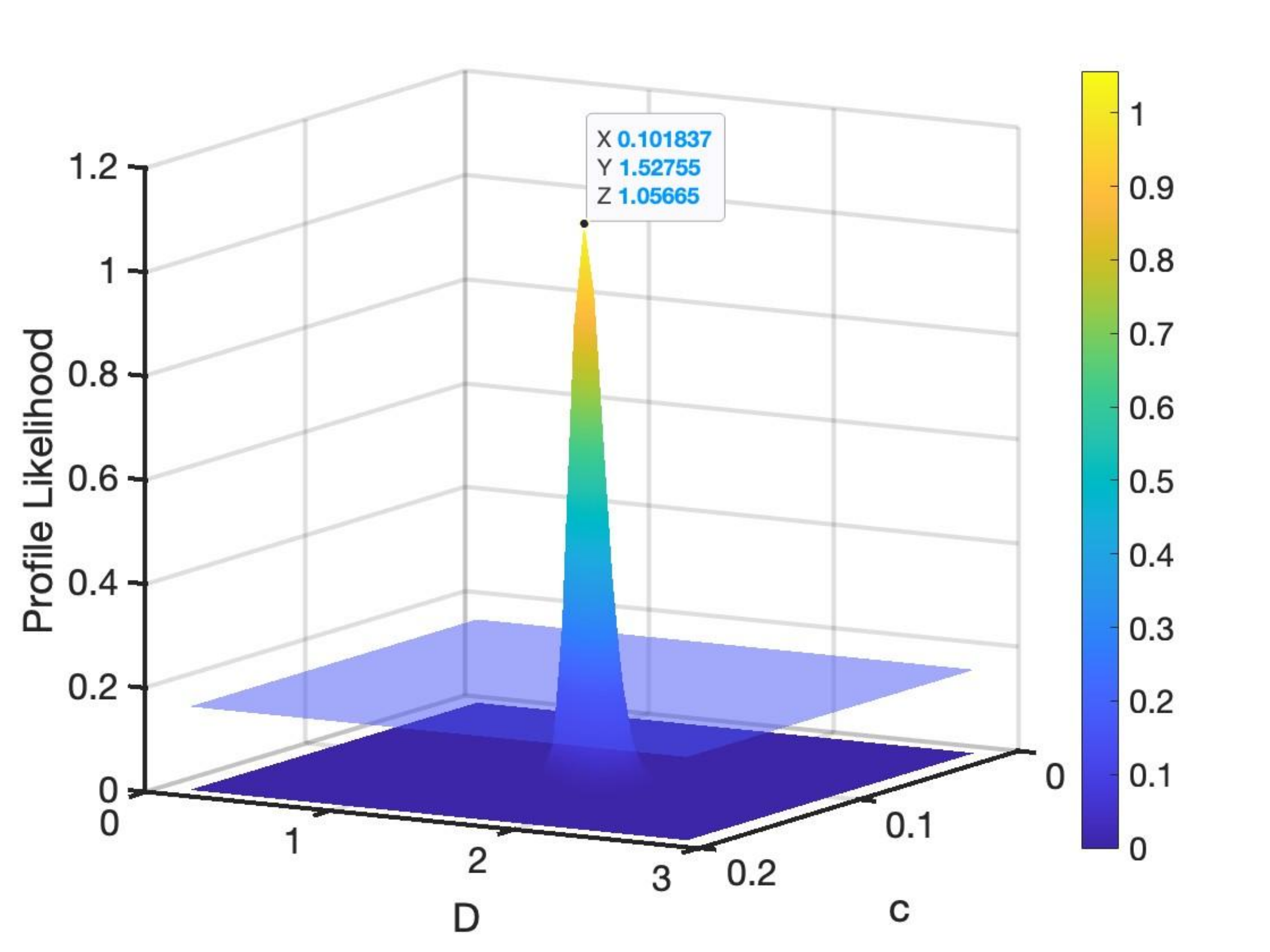}
    \end{minipage}
    
    \vspace{0.5em}
    \begin{minipage}[t]{0.485\linewidth}
        \centering
        \includegraphics[width=\linewidth]{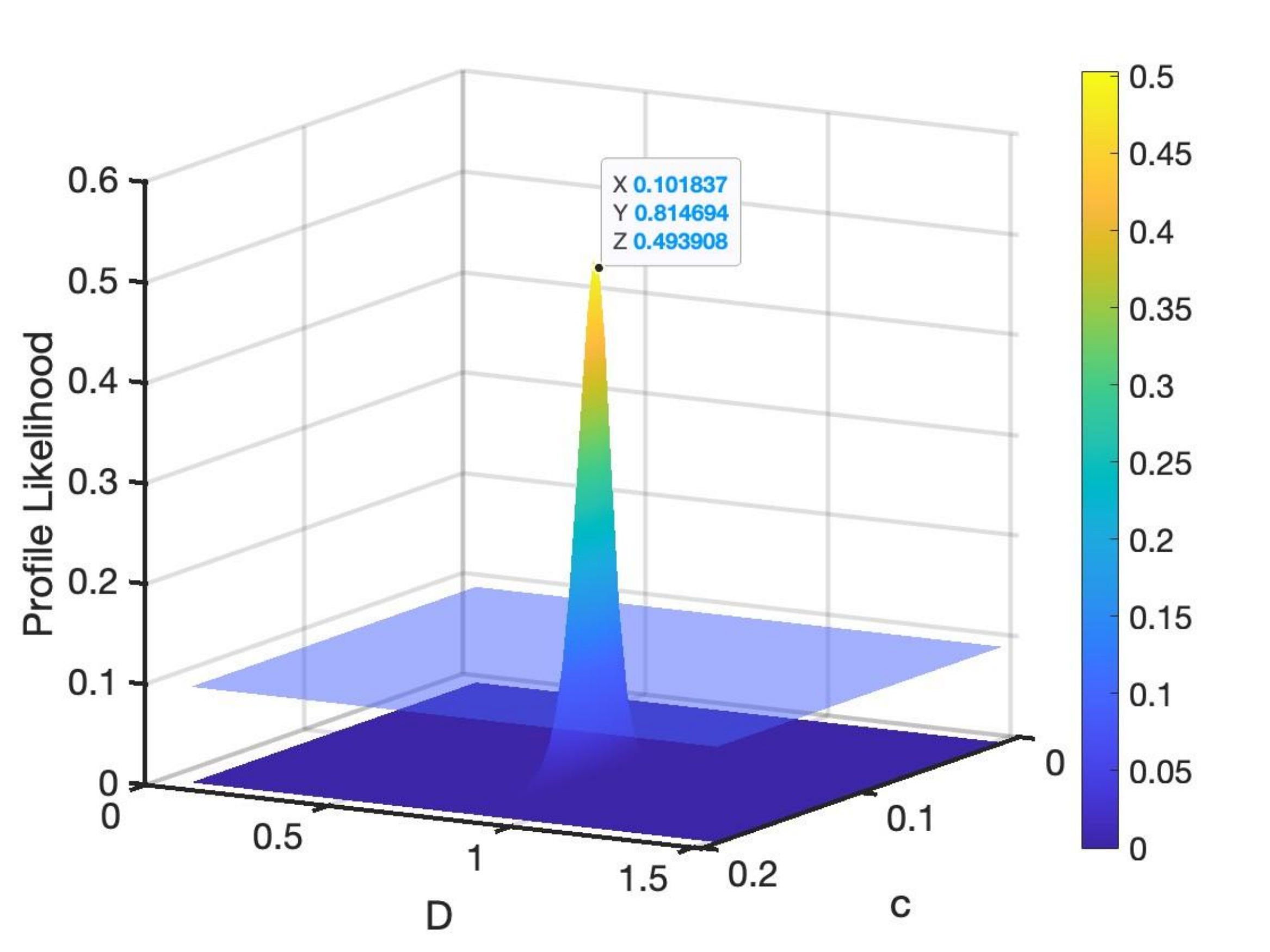}
    \end{minipage}
    
    \caption{Two-dimensional profile likelihood given each possible pair of $(c, D)$ chosen around their baseline values in each region. The z-axis shows the corresponding profile likelihood values calculated using \Cref{eq: PL}, together with a threshold plane based on \Cref{eq: CI3} (top-left: Region~1, top-right: Region~2, bottom: Region~3).}
    \label{fig:cD_PL}
\end{figure}

We visualize the two-dimensional profile likelihood results for the three regions in our model and data by drawing the corresponding 3D plots in \Cref{fig:cD_PL}. For each point on the $(c, D)$-plane, we obtain one value for the least square error, which is computed by calculating the difference between the synthetic data curve generated under optimized values of $\beta_1$ and $\beta_2$ and the baseline value curve (red curves shown in \Cref{tab:baseline}). The likelihood value is then derived by \Cref{eq: PL}. It can be seen from \Cref{fig:cD_PL} that the profile likelihood value achieves its peak around the baseline values for $c$ and $D$ as shown in \Cref{tab:baseline}. A transparent plane is added in each plot, which represents the threshold calculated using \Cref{eq: CI3}. 

Overall, we find that both $c$ and $D$ are identifiable, which is consistent with the one-dimensional profile likelihood results. Two-dimensional profile likelihood thus offers a valuable tool for analyzing experimental data by providing both the range of plausible measurement outcomes for a given experiment and their impact on the parameter likelihood profile \cite{litwin2022optimal,kreutz2013profile}. However, a limitation of this approach is that it is computationally more intensive compared to one-dimensional analyses, especially when applied to models with many parameters. Nevertheless, our results demonstrate that for the FRAP data and cellular regions considered here, the transport speed $c$ and diffusion coefficient $D$ can be assumed to be identifiable.


\section{Investigating parameter relationships in transport models of FRAP experiments}
\subsection{Subset profiles and 3D contour plots}
According to the results derived from the one-dimensional profile likelihood, we find that $c$ and $D$ are identifiable, while $\beta_1$ and $\beta_2$ are non-identifiable. Then, we expand the methodology in \cite{ciocanel2024parameter} in order to investigate potential parameter relationships between $\beta_1$ and $\beta_2$ that may be identifiable.

\begin{figure}[htbp]
    \centering
    \begin{minipage}[t]{0.485\linewidth}
        \centering
        \includegraphics[width=\linewidth]{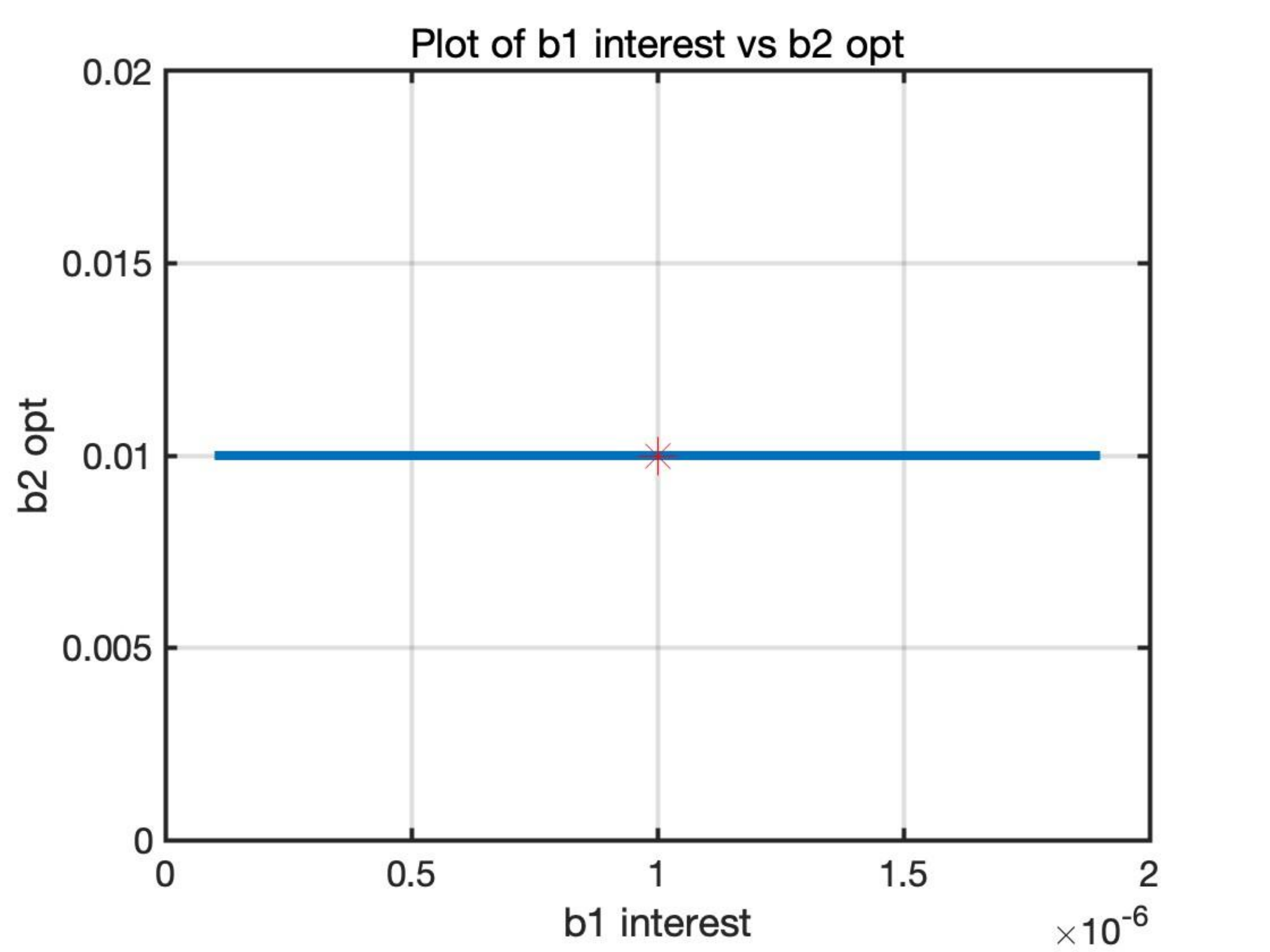}
    \end{minipage}
    \hspace{0.01\linewidth}
    \begin{minipage}[t]{0.485\linewidth}
        \centering
        \includegraphics[width=\linewidth]{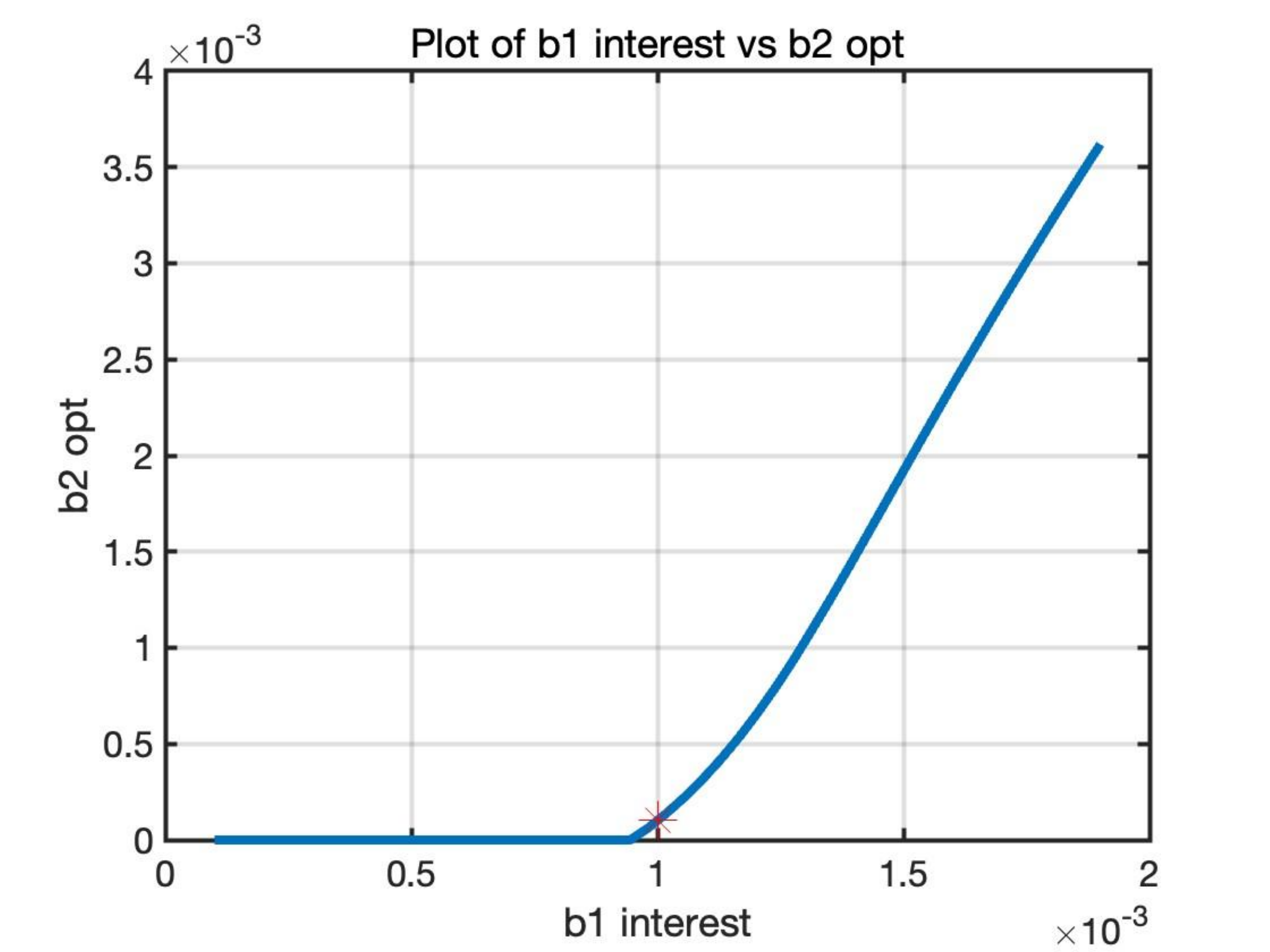}
    \end{minipage}
    
    \vspace{0.5em}
    \begin{minipage}[t]{0.485\linewidth}
        \centering
        \includegraphics[width=\linewidth]{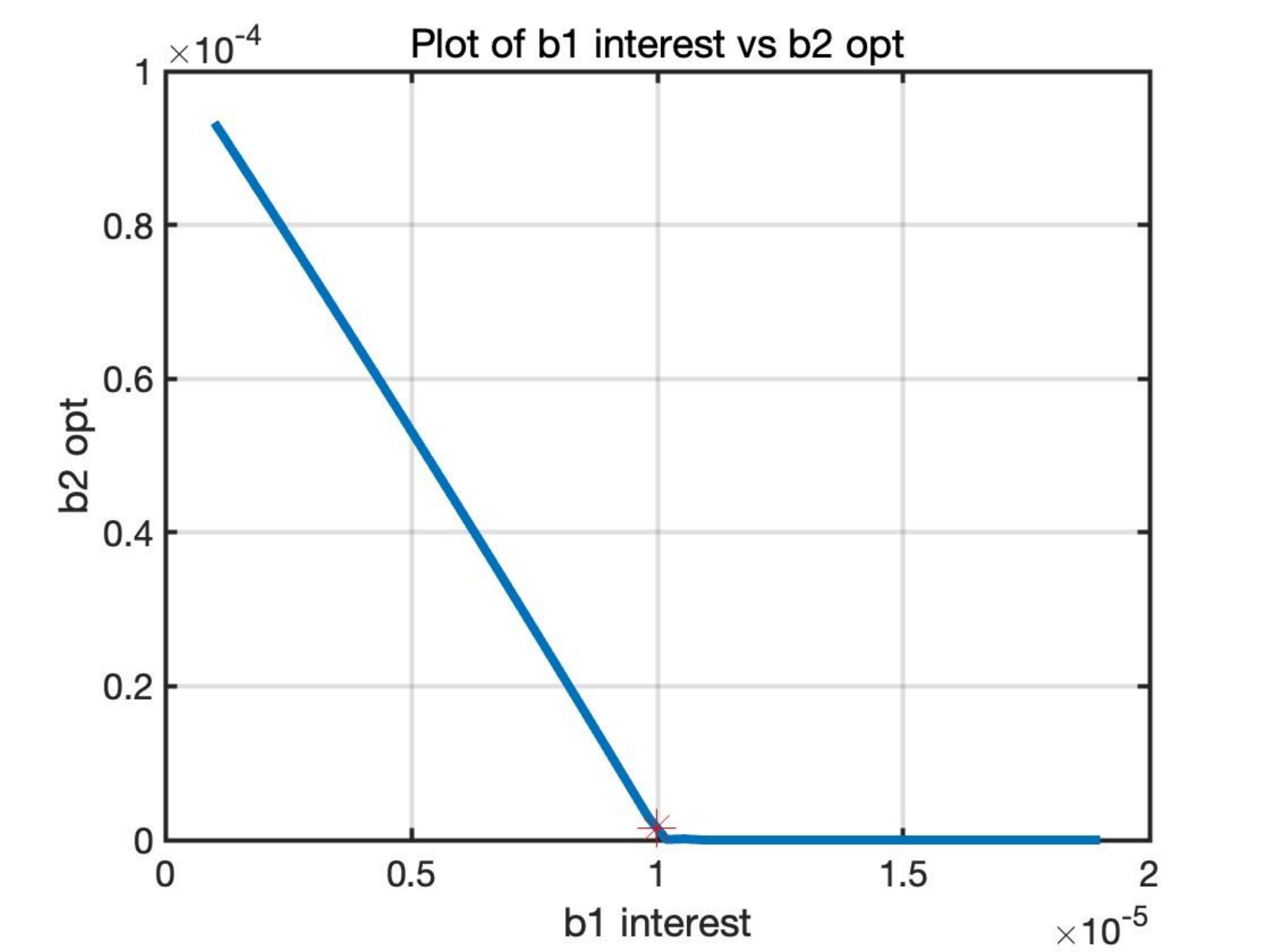}
    \end{minipage}
    
    \caption{Subset profiles for interest parameter $\beta_1$ on the x-axis and the corresponding optimized $\beta_2$ on the y-axis for each region given noiseless FRAP data synthetically generated using \Cref{eq: PDEs} and the parameter set in  \Cref{tab:baseline}. Baseline values are indicated as red stars in each plot (top-left: Region~1, top-right: Region~2, bottom: Region~3).}
    \label{fig:subset}
\end{figure}

We first begin by drawing the subset profiles based on the data from profile likelihood for the three regions. In each subset profile, we treat $\beta_1$ as the interest parameter. For each fixed value of $\beta_1$, we record the corresponding optimal value of $\beta_2$ that maximize the likelihood value. This allows us to visualize the compensatory relationship between these two kinetic parameters \cite{ciocanel2024parameter}. As shown in \Cref{fig:subset}, the x-axis refers to the values for interest parameter $\beta_1$, and the y-axis denotes the corresponding optimal values for $\beta_2$.


\begin{figure}[htbp]
    \centering
    \begin{minipage}[t]{0.485\linewidth}
        \centering
        \includegraphics[width=\linewidth]{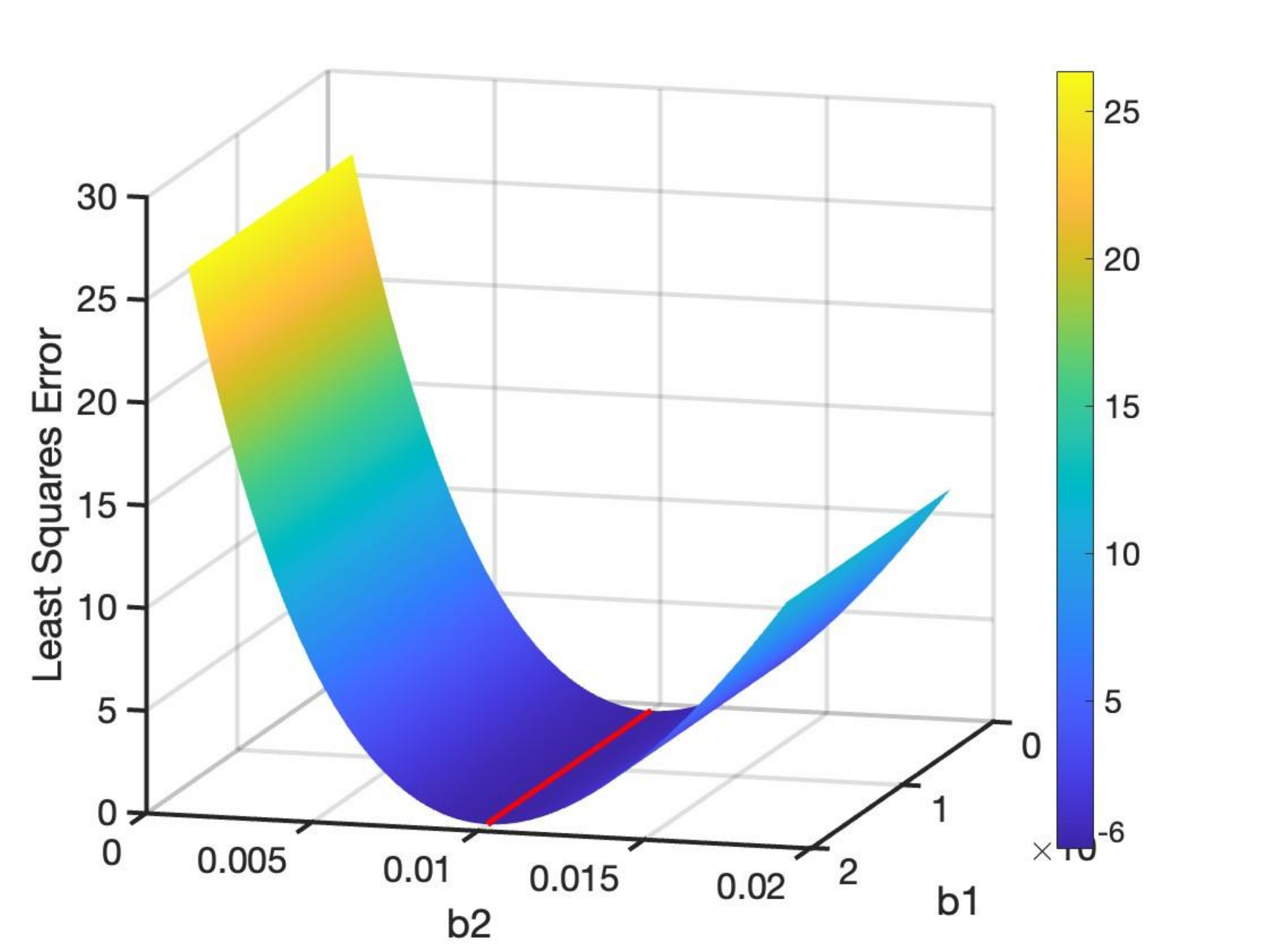}
    \end{minipage}
    \hspace{0.01\linewidth}
    \begin{minipage}[t]{0.485\linewidth}
        \centering
        \includegraphics[width=\linewidth]{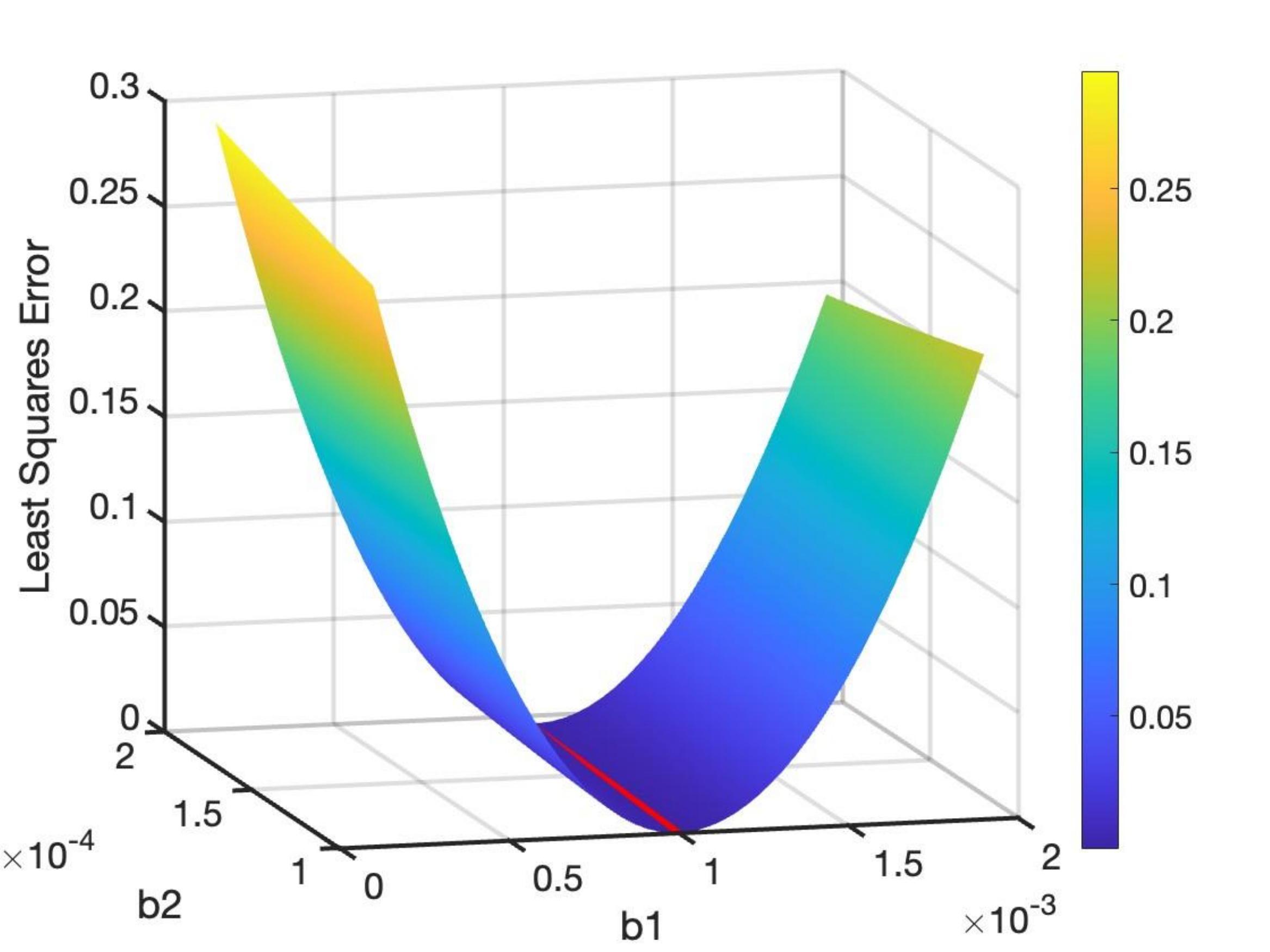}
    \end{minipage}
    
    \vspace{0.5em}
    \begin{minipage}[t]{0.485\linewidth}
        \centering
        \includegraphics[width=\linewidth]{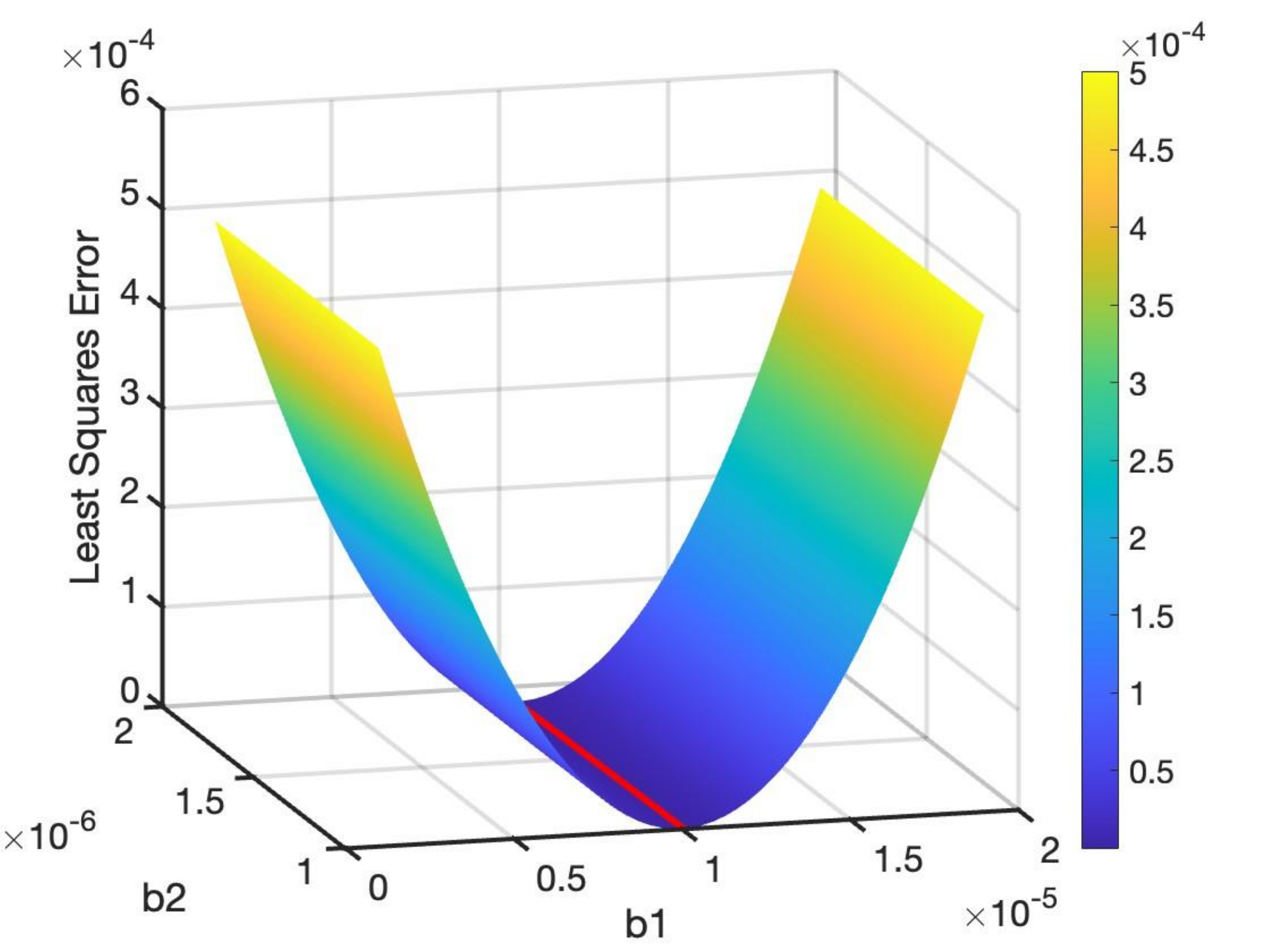}
    \end{minipage}
    
    \caption{Contour plots of the least square error between synthetic FRAP data generated using baseline values of $c$ and $D$ in \Cref{tab:baseline}, as well as $\beta_1$ and $\beta_2$ values from a grid around their baseline values for each region (top-left: Region~1, top-right: Region~2, bottom: Region~3).}
    \label{fig:3d}
\end{figure}

\Cref{fig:subset} shows that there could be some linear relationships between $\beta_1$ and $\beta_2$ in some regions of parameter space. To better interpret these relationships, our goal is to explore the likelihood landscape by varying $\beta_1$ and $\beta_2$ in a grid around their baseline values in each region under fixed values for $c$ and $D$ (set to their baseline values in \Cref{tab:baseline}). We then carry out the forward PDE computation by plugging the four chosen parameter values into PDE model (\Cref{eq: PDEs}) to obtain another synthetic data curve. For each point in the $\beta_1$-$\beta_2$ plane, there is one corresponding least square error by calculating the difference between the synthetic data curve and the baseline curve (red curves shown in \Cref{fig:combined1}). \Cref{fig:3d} shows the contour plot of this error as a function of ($\beta_1$, $\beta_2$) for the three regions. We notice that the least square error is minimized along a line (highlighted in red). Therefore, we cannot distinguish parameters $(\beta_1, \beta_2)$ on this red line, as they all give the minimized error function \cite{ciocanel2024parameter}.


\subsection{Slope vector field generation}
In order to find the corresponding relationships between $\beta_1$ and $\beta_2$ in each region, we generate a slope vector field as in \cite{ciocanel2024parameter}. In the slope vector field, each vector represents the slope in the subset profiles showing the relationship between $\beta_1$ and $\beta_2$ (\Cref{fig:subset}). The goal is to draw a contour curve in the vector field for each region, representing the combination of $\beta_1$ and $\beta_2$ that yields the maximum profile likelihood value.

\begin{figure}[htbp]
    \centering
    \begin{minipage}[t]{0.485\linewidth}
        \centering
        \includegraphics[width=\linewidth]{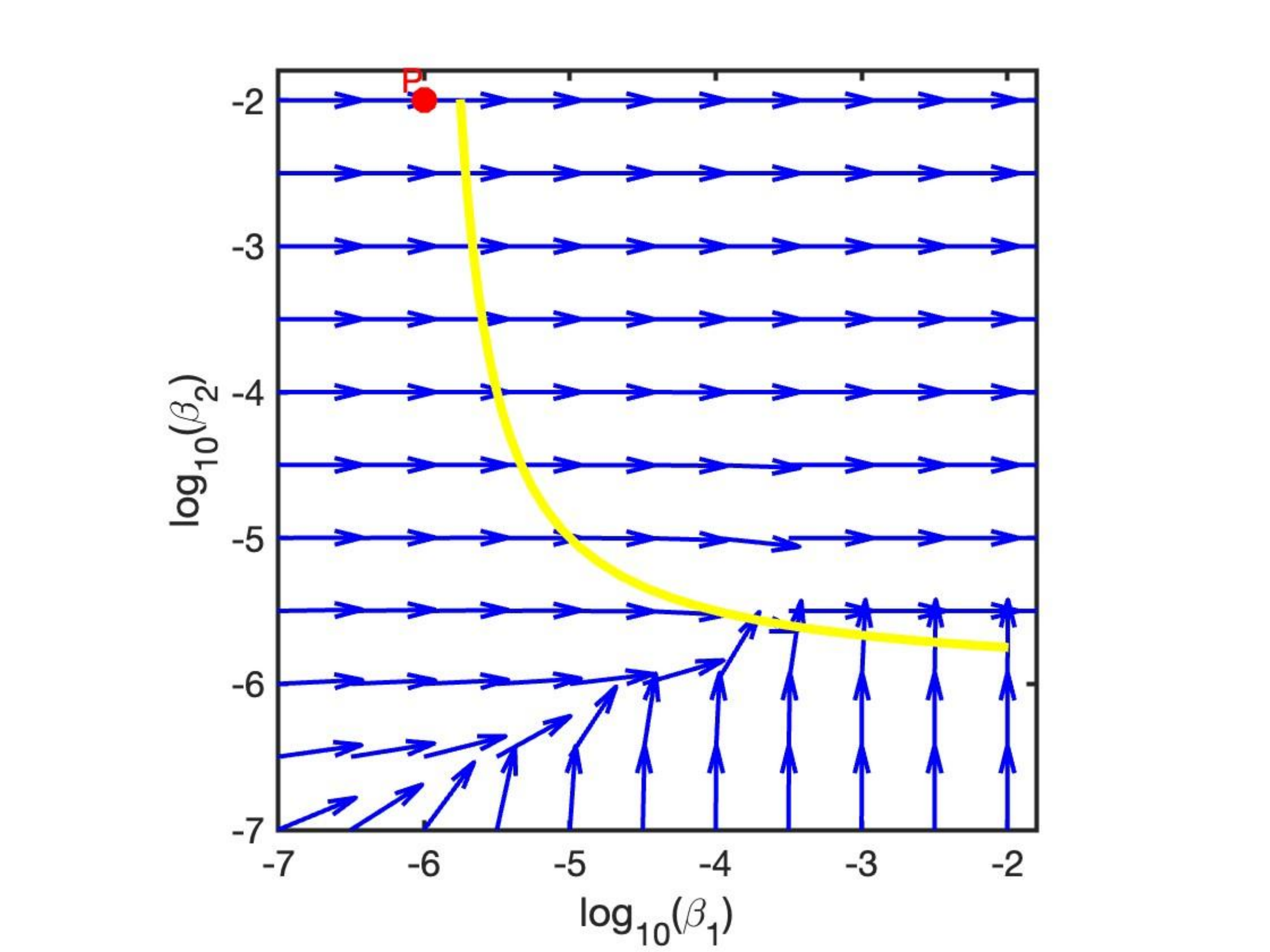}
    \end{minipage}
    \hspace{0.01\linewidth}
    \begin{minipage}[t]{0.485\linewidth}
        \centering
        \includegraphics[width=\linewidth]{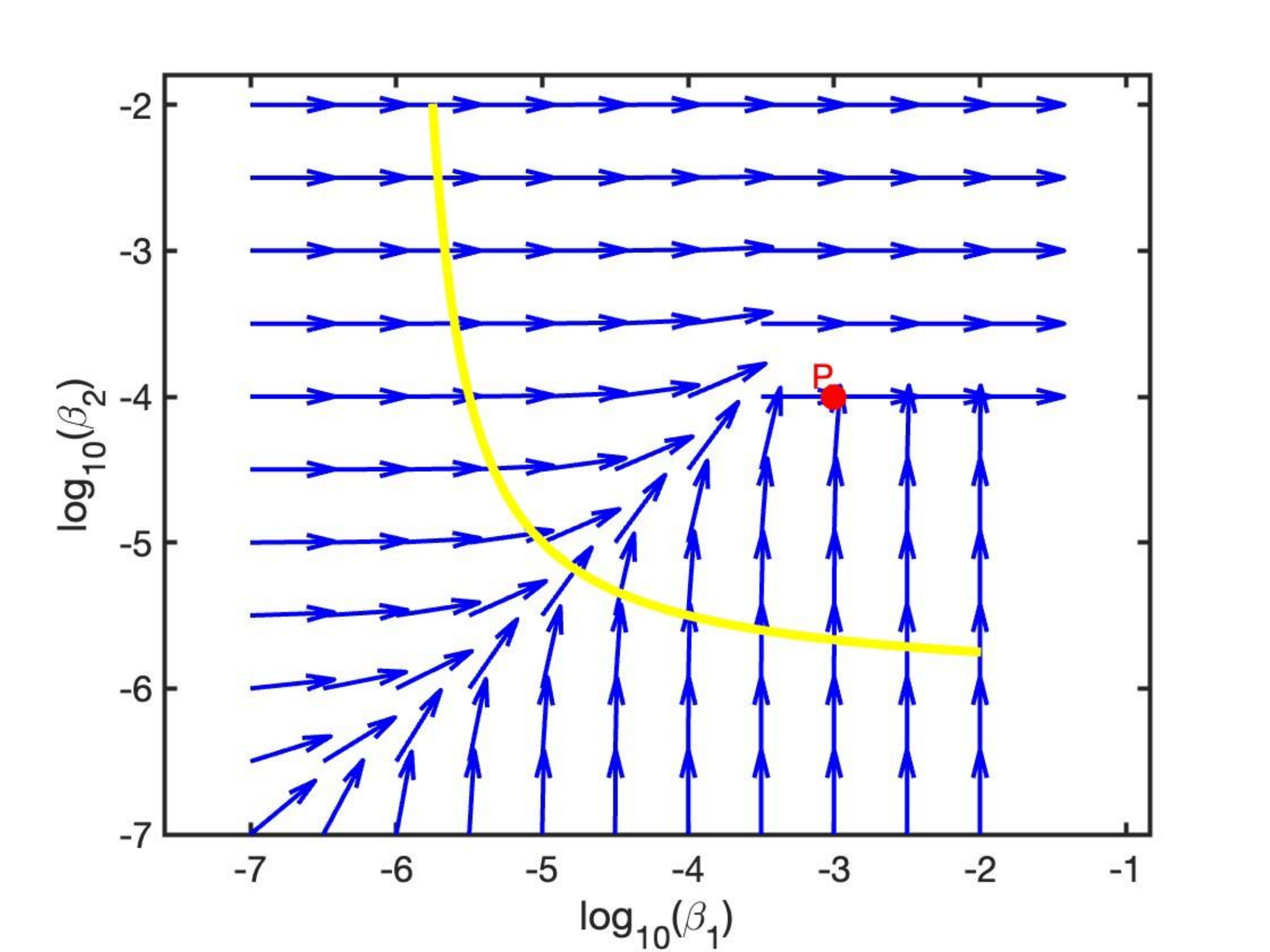}
    \end{minipage}
    
    \vspace{0.5em}
    \begin{minipage}[t]{0.485\linewidth}
        \centering
        \includegraphics[width=\linewidth]{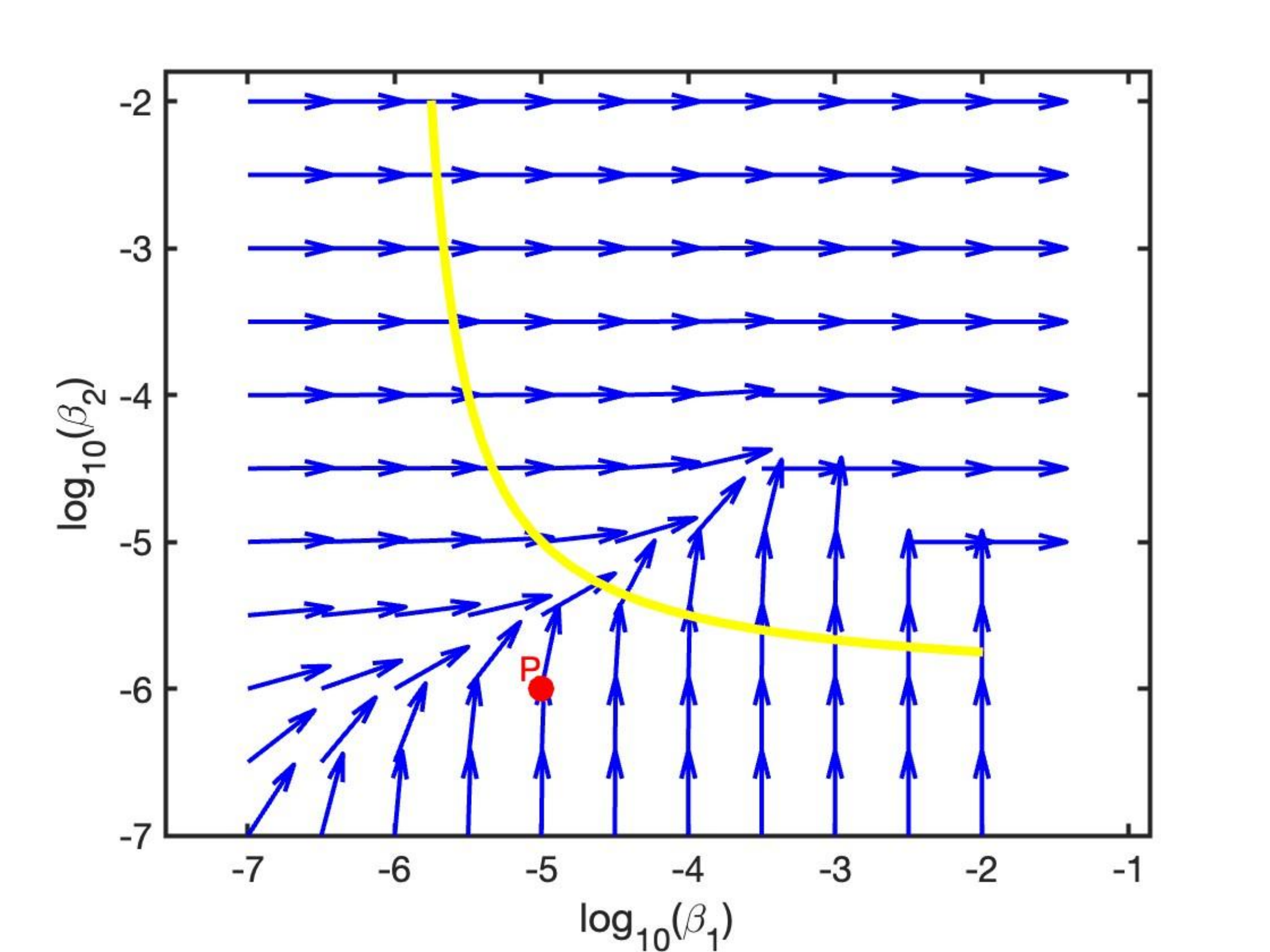}
    \end{minipage}
    
    \caption{Slope vector field generation: the slope of each arrow corresponds to the slope in the subset profile for each profile likelihood. A yellow curve $\tau$ crosses transversely the slope vector field, with point $P$ indicating the baseline values for $\beta_1$ and $\beta_2$ (top-left: Region~1, top-right: Region~2, bottom: Region~3).}
    \label{fig:vector}
\end{figure}

Since the diffusion coefficient $D$ and transport speed $c$ are identifiable, we set them to be their baseline values in \Cref{tab:baseline}. Then, we select a grid in ($log_{10}{\beta_1}$, $log_{10}{\beta_2}$)-plane. For each point in the plane, we generate a synthetic FRAP data set based on \Cref{eq: PDEs} and the chosen parameter values. For each synthetic data set, we compute the profile likelihood for $\beta_1$ and $\beta_2$, then set the tangent vector at each point to be the slope of the subset profile shown in \Cref{fig:subset}. The pattern of the slope vector field for each region is shown in \Cref{fig:vector}.

We then choose a curve $\tau$ going transversely across the slope vector field, meaning that for each point on this curve, its tangent vector is not aligned with the direction of the vector field. We can either choose the curve $\tau$ in explicit analytical form or else using linear interpolation and a forward Euler scheme applied to gradients of the vector to compute such a curve numerically (\cite{ciocanel2024parameter}). Here, we use an explicit analytical parametrization similar to the one in \cite{ciocanel2024parameter}:

\begin{equation}
\begin{aligned}
    \log_{10}{\beta_1} &= s + \sqrt{s^2+1}-6 \\
    \log_{10}{\beta_2} &= -s + \sqrt{s^2+1}-6
\end{aligned}
\label{eq:beta12}
\end{equation}

This reparametrization yields the yellow curve for each region shown in \Cref{fig:vector}. Moreover, in each plot, the point $P$ corresponds to the baseline values of $\beta_1$ and $\beta_2$ in each region.

\begin{figure}[htbp]
    \centering
    \begin{minipage}[t]{0.485\linewidth}
        \centering
        \includegraphics[width=\linewidth]{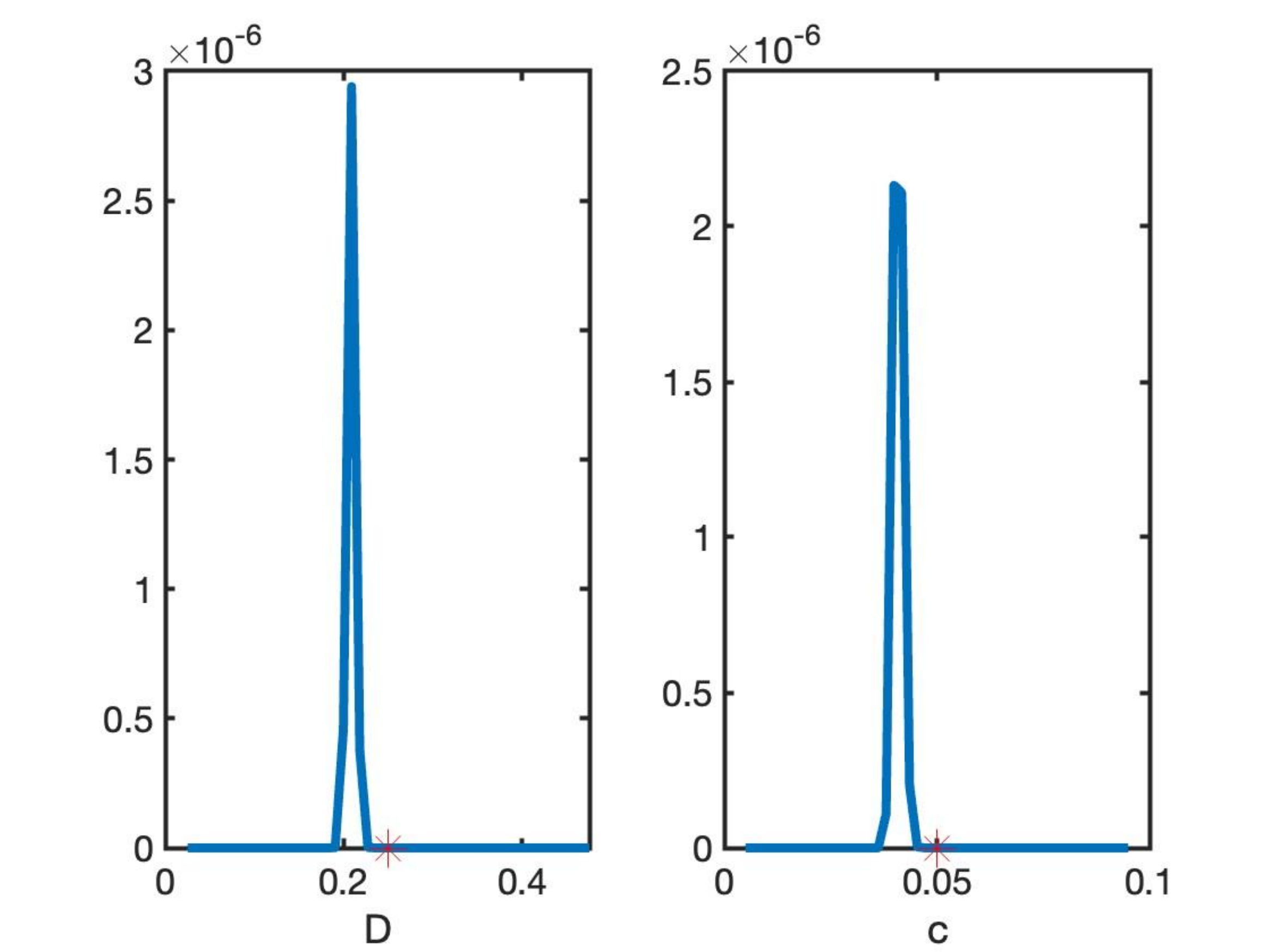}
    \end{minipage}
    \hspace{0.01\linewidth}
    \begin{minipage}[t]{0.485\linewidth}
        \centering
        \includegraphics[width=\linewidth]{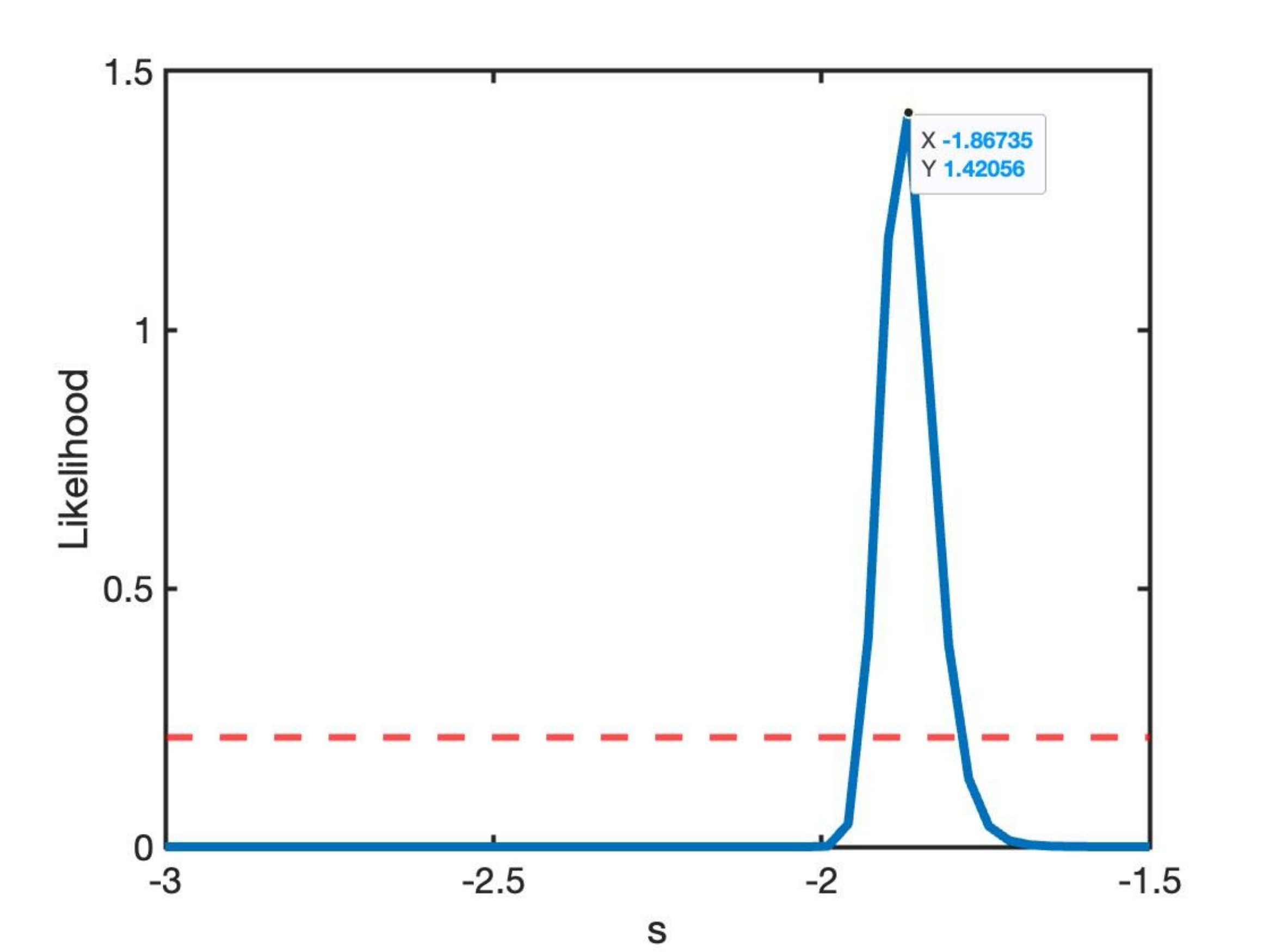}
    \end{minipage}
    
    \vspace{0.5em}
    \begin{minipage}[t]{0.7\linewidth}
        \centering
        \includegraphics[width=\linewidth]{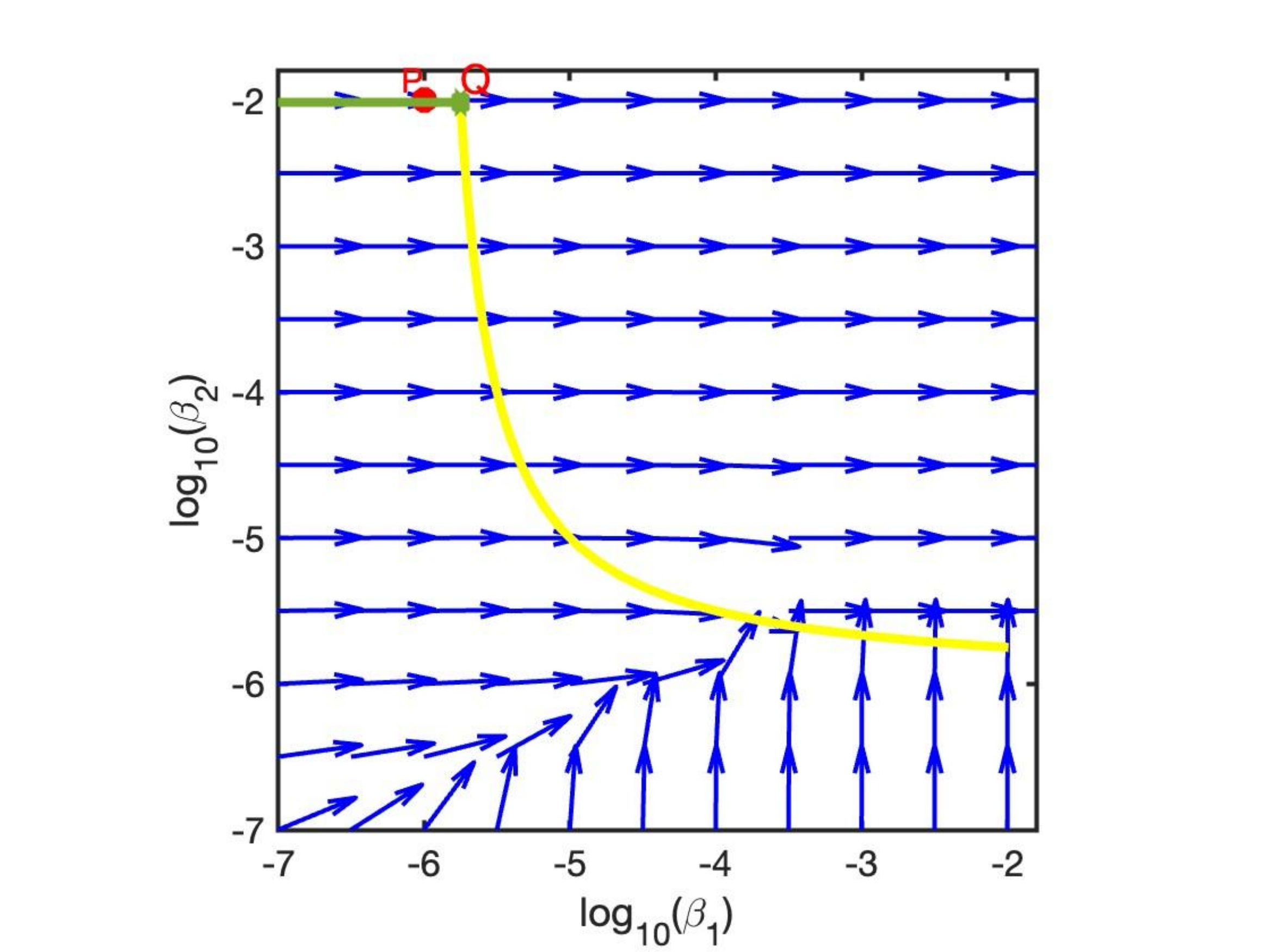}
    \end{minipage}
    
    \caption{Profile likelihoods for each interest parameter $c$, $D$, and $s$ in Region~1 given noiseless FRAP data generated using the parameters in \Cref{tab:baseline} (top-left and top-right). $s$ achieves the maximum at $s^* = -1.86735$, corresponding to the point $Q$ on the yellow curve $\tau$, as well as the trace of the error-minimizing green contour curve (bottom).}
    \label{fig:r1}
\end{figure}

\begin{figure}[htbp]
    \centering
    \begin{minipage}[t]{0.485\linewidth}
        \centering
        \includegraphics[width=\linewidth]{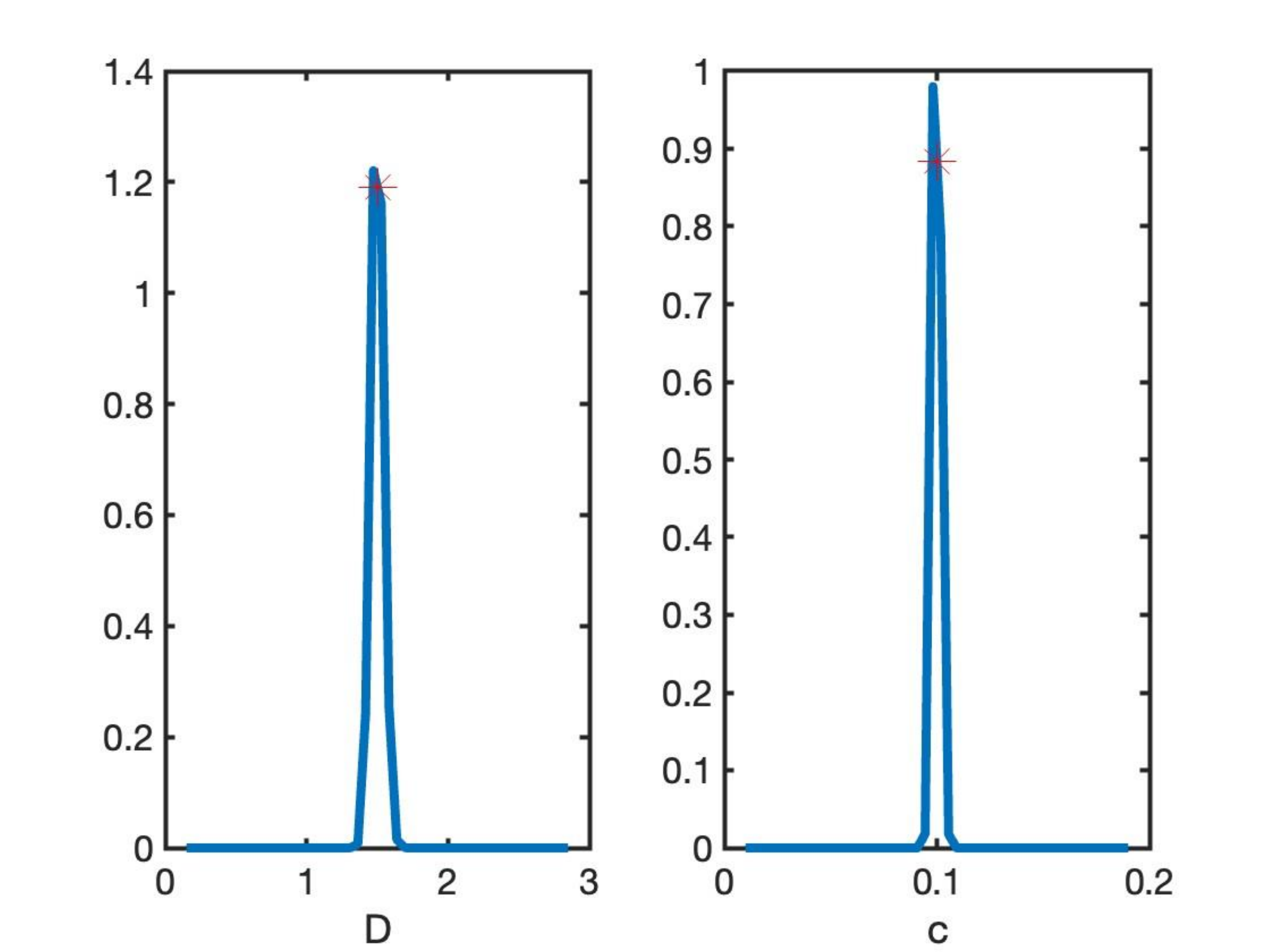}
    \end{minipage}
    \hspace{0.01\linewidth}
    \begin{minipage}[t]{0.485\linewidth}
        \centering
        \includegraphics[width=\linewidth]{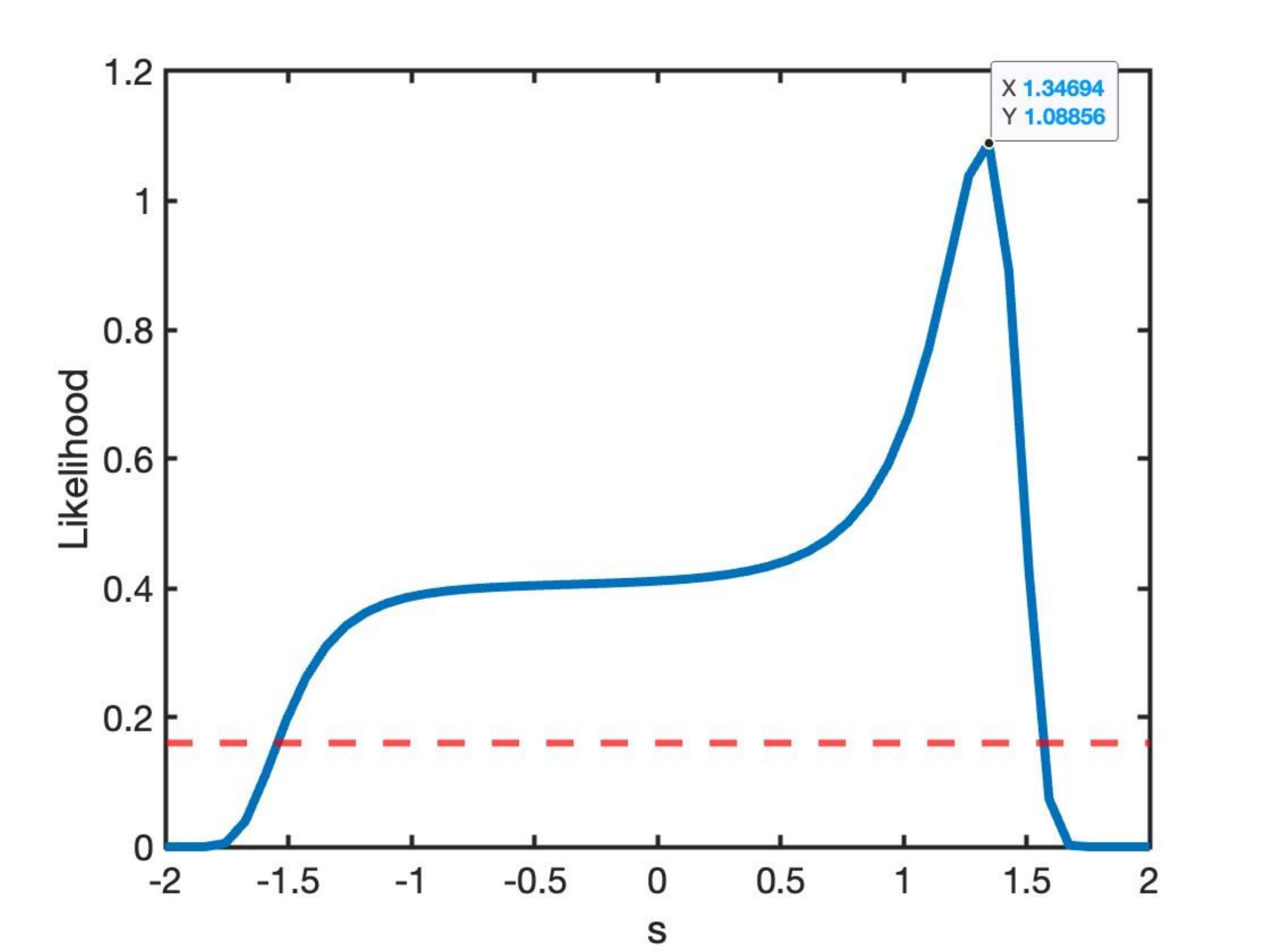}
    \end{minipage}

    \vspace{0.5em}
    \begin{minipage}[t]{0.7\linewidth}
        \centering
        \includegraphics[width=\linewidth]{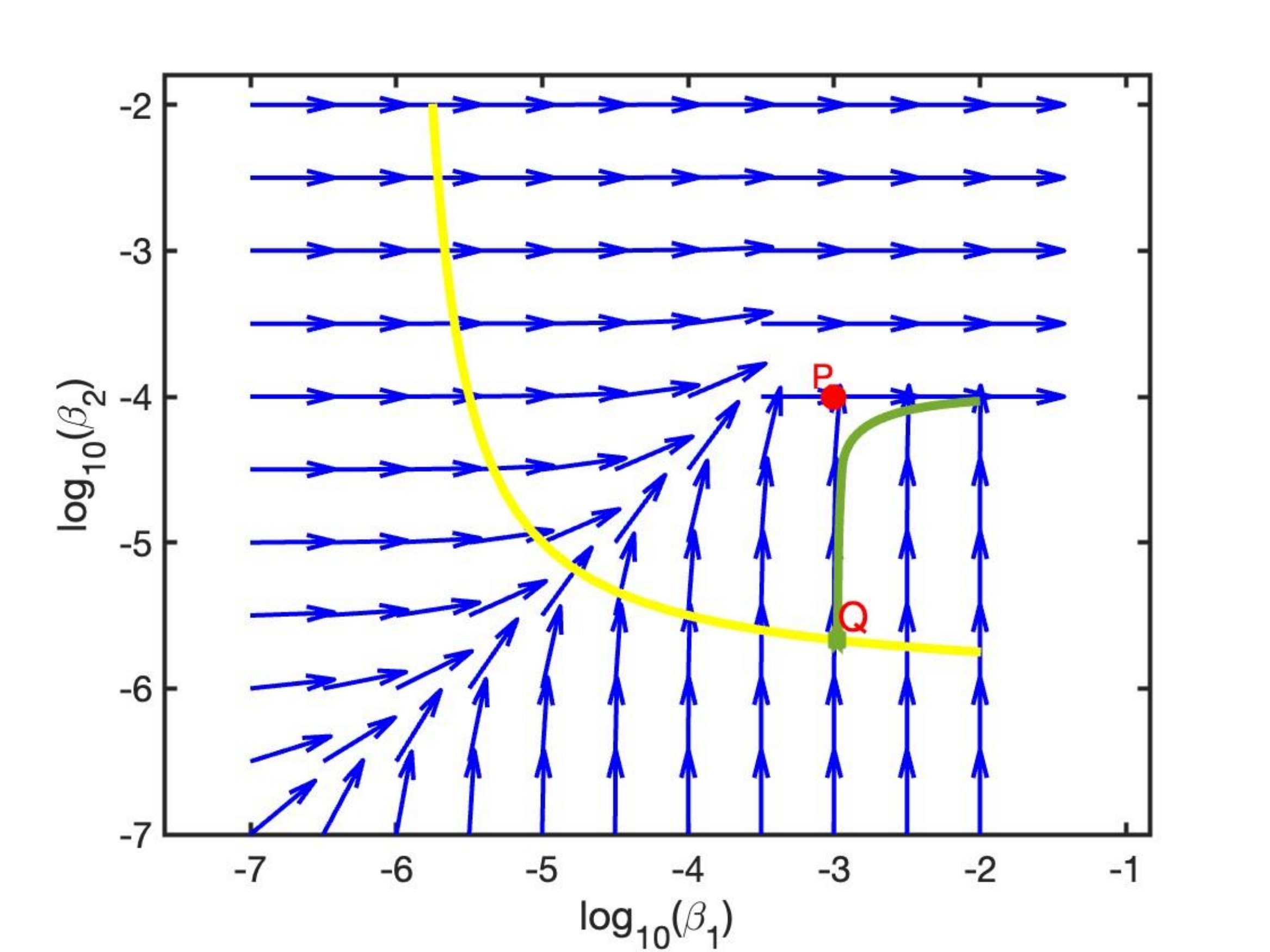}
    \end{minipage}

    \caption{Profile likelihoods for each interest parameter $c$, $D$ and $s$ in Region 2 given noiseless FRAP data generated using the parameters in \Cref{tab:baseline} (top-left and top-right). $s$ achieves the maximum at $s^* = 1.34694$, corresponding to the point $Q$ on the yellow curve $\tau$, as well as the trace of error-minimizing green contour curve (bottom).}
    \label{fig:r2}
\end{figure}

Then, we apply the profile likelihood analysis for the three parameters $(c, D,s)$. We especially want to see whether the new parameter $s$ is identifiable or not. First, both $c$ and $D$ are identifiable, which is consistent with the results shown in \Cref{fig:PL_all}. As expected, in \Cref{fig:r1} and \Cref{fig:r2}, the clear peaks of the s-profile likelihood for Region 1 and Region 2 demonstrate that the re-parametrized parameter $s$ is indeed identifiable. Next, we plug the optimized $s$ value at the peak into \Cref{eq:beta12} to get the values for $\beta_1$ and $\beta_2$, and that will correspond to a point $Q$ on the curve $\tau$. By the method of linear interpolation, we then generate a green contour curve which begins from $Q$. For each point on this curve, the profile likelihood function is maximized. Therefore, this green contour curve represents the relationship between $\beta_1$ and $\beta_2$ that is identifiable in both regions. Notably, the ground-truth point $p$ is very close to this contour curve.

\begin{figure}[htbp]
    \centering
    \begin{minipage}[t]{0.485\linewidth}
        \centering
        \includegraphics[width=\linewidth]{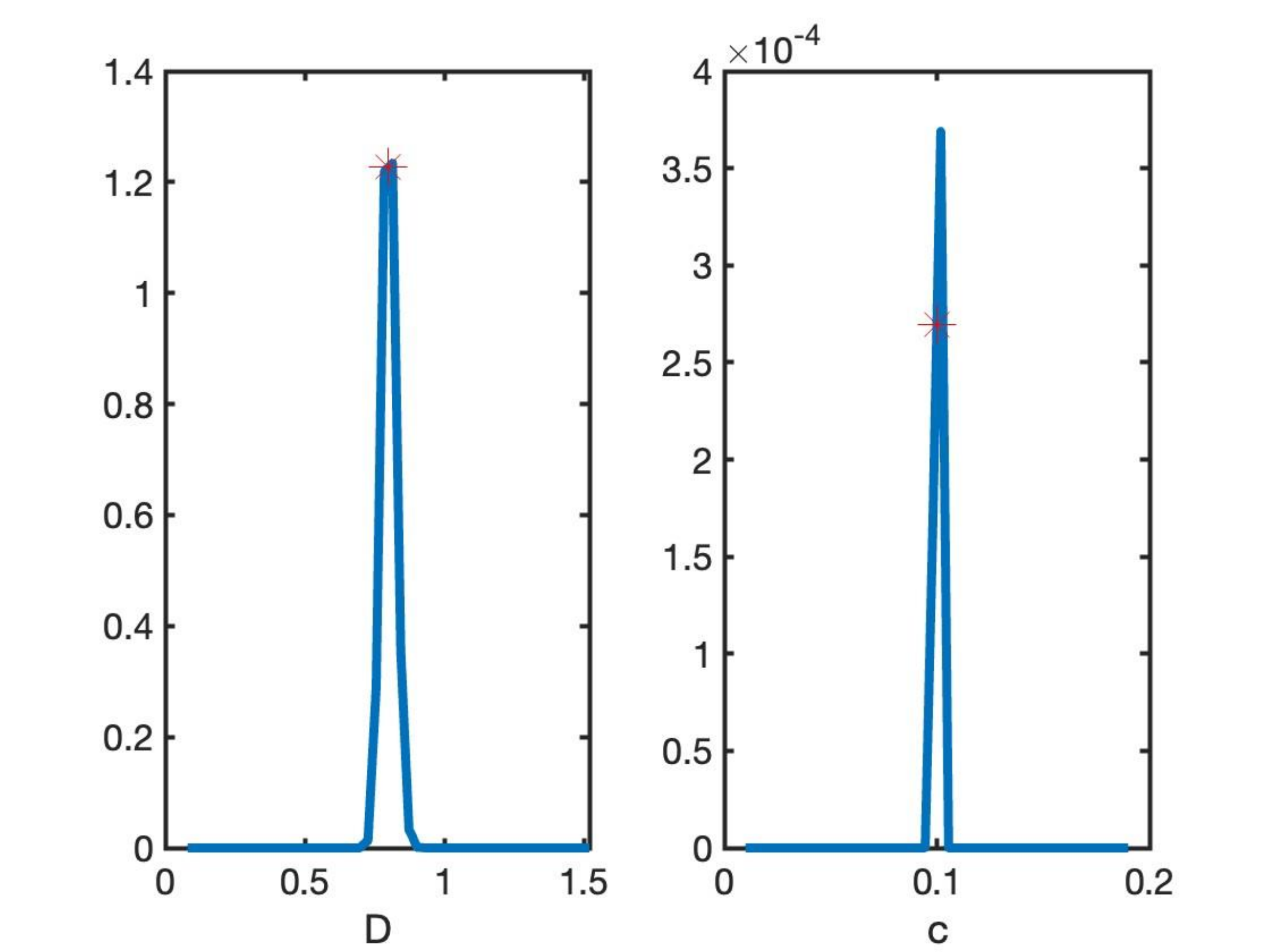}
    \end{minipage}
    \hspace{0.01\linewidth}
    \begin{minipage}[t]{0.485\linewidth}
        \centering
        \includegraphics[width=\linewidth]{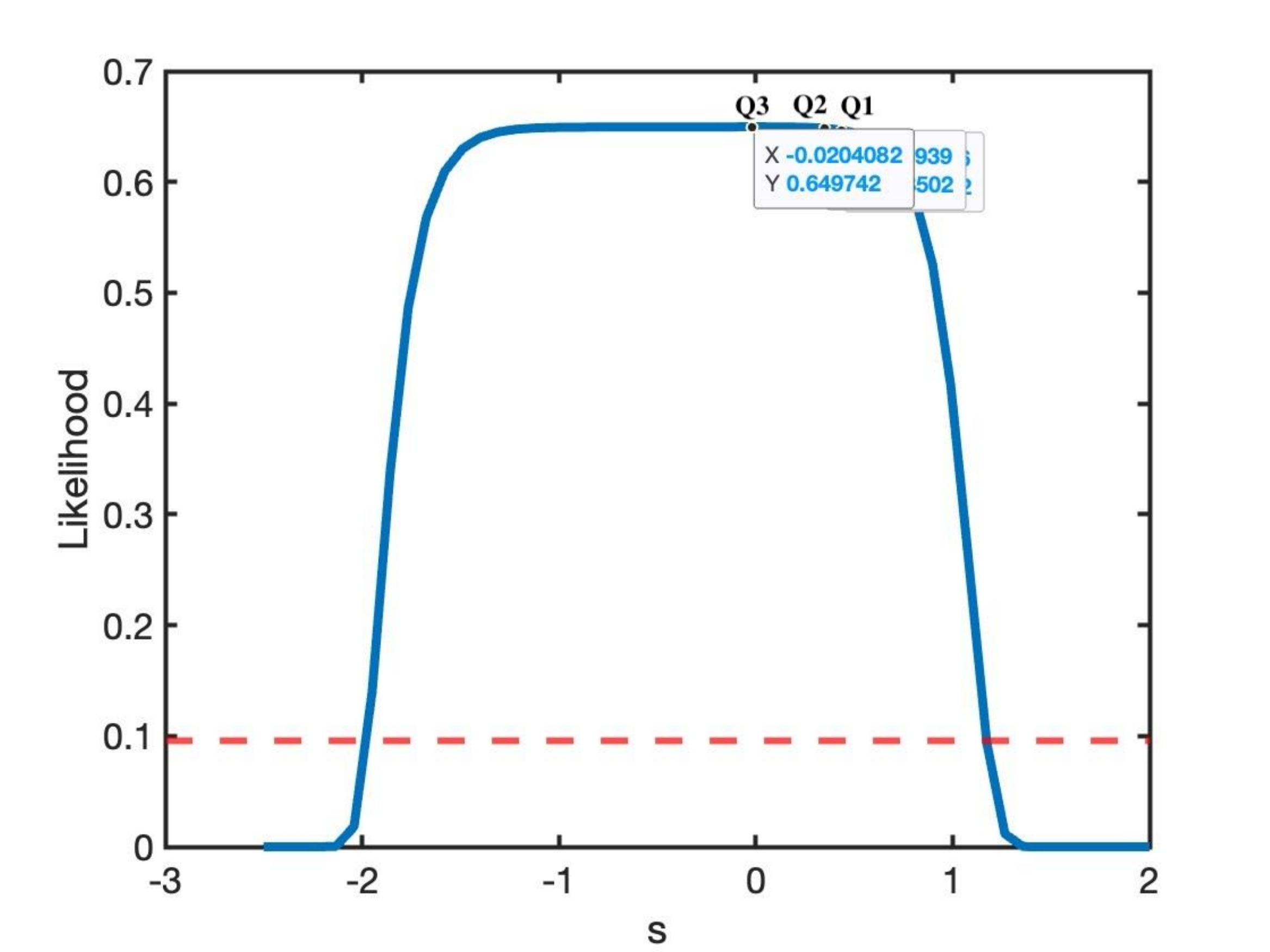}
    \end{minipage}

    \vspace{0.5em}
    \begin{minipage}[t]{0.7\linewidth}
        \centering
        \includegraphics[width=\linewidth]{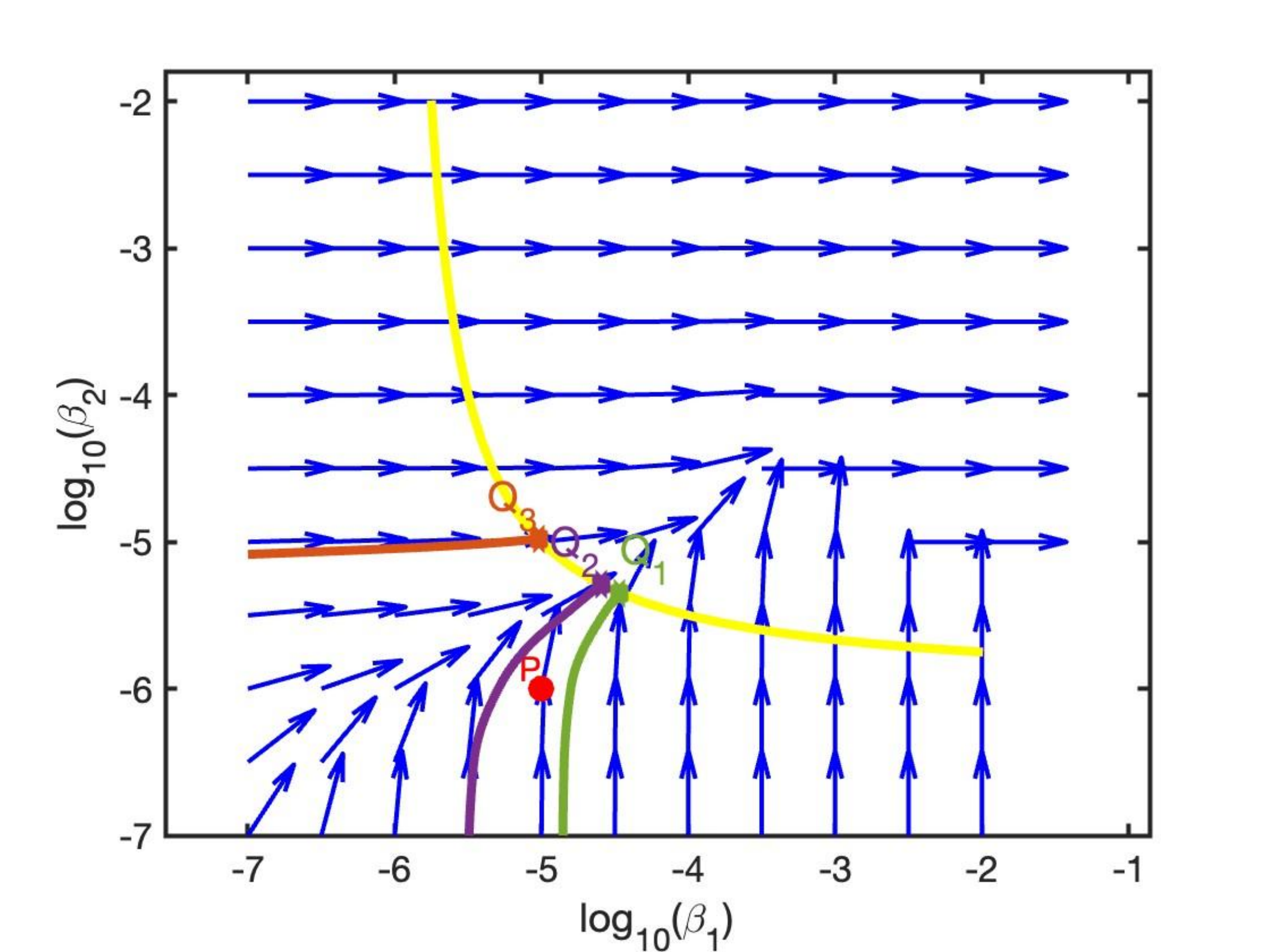}
    \end{minipage}

    \caption{Profile likelihoods for each interest parameter $c$, $D$, and $s$ in Region~3 given noiseless FRAP data generated using the parameters in \Cref{tab:baseline} (top-left and top-right). Choosing three values for $s$: $s^* = 0.438776, 0.346939, -0.0204082$ corresponds to points $Q_1, Q_2, Q_3$ on the yellow curve $\tau$, as well as the traces of error-minimizing contour curves (bottom).}
    \label{fig:r3}
\end{figure}

For Region 3, the likelihood value is maximized along a line instead of one peak as shown in \Cref{fig:r3}. Therefore, we plug in some $s$ values on that line into the reparametrization formula, and that corresponds to points $Q_1$, $Q_2$, $Q_3$ on the curve $\tau$. Based on each $Q$, we then generate the corresponding contour curves. Since we know that the ground-truth point $P$ must be very close to the contour curve, the green and the purple curve beginning from $Q_1$ and $Q_2$ are more likely to bound the relationship between $\beta_1$ and $\beta_2$. However, we cannot determine a contour curve that can best represent the relationship between $\beta_1$ and $\beta_2$ in this case, since there is not a single peak as depicted in the s-profile likelihood plot. This means that, for region 3, we are not able to determine an identifiable correlation for $\beta_1$ and $\beta_2$ without additional information.


\subsection{Pipeline for identifiability analysis and parameter relationship investigation}
\label{subsec:pipeline}

Here, we outline the systematic workflow employed to analyze parameter identifiability and explore relationships among non-identifiable parameters. The pipeline consists of four key steps that progressively refine our understanding of the model parameters. We will apply this
pipeline to real experimental data in \Cref{application}.

\begin{enumerate}

\item \textbf{Profile likelihood computation}

The first step assesses identifiability of individual parameters using synthetic FRAP data. For each parameter of interest (e.g., $c$, $D$, $\beta_1$, $\beta_2$), we fix its value across a grid around the baseline values (\Cref{tab:baseline}) while optimizing the remaining nuisance parameters to maximize the profile likelihood (\Cref{eq: PL}). We establish the 95\% confidence threshold (\Cref{eq: CI3}) to judge identifiability. This reveals that identifiable parameters like $c$ and $D$ exhibit sharp peaks in their likelihood profiles, while non-identifiable parameters (e.g., $\beta_1$, $\beta_2$) show flat or plateaued profiles (\Cref{fig:PL_all}).

\item \textbf{Two-dimensional profile likelihood exploration}

The second step validates pairwise identifiability of parameters $c$ and $D$. We define a grid of $(c, D)$ values around baseline estimates (\Cref{tab:baseline}) and for each pair, optimize $\beta_1$ and $\beta_2$ to maximize likelihood (\Cref{eq: 2D PL}). The resulting likelihood landscape visualized as a 3D surface (\Cref{fig:cD_PL}) shows a clear peak near baseline values, confirming that $c$ and $D$ are jointly identifiable.

\item \textbf{Contour interpretation for non-identifiable parameters}

The third step investigates relationships between $\beta_1$ and $\beta_2$ when $c$ and $D$ are fixed at baseline values. We compute the least square error (LSE) for a $(\beta_1, \beta_2)$ grid and identify the contour line where the LSE is minimized (\Cref{fig:3d}). This reveals a correlation between $\beta_1$ and $\beta_2$ along the minimized LSE line, indicating structural non-identifiability due to compensatory effects between these parameters.

\item \textbf{Slope vector field and reparametrization}

The final step derives identifiable combinations of $\beta_1$ and $\beta_2$:

\begin{itemize}
    \item For subset profile construction, for each $\beta_1$ value on a logarithmic grid, we compute the optimal $\beta_2$ that maximizes the profile likelihood while fixing $c$ and $D$ at baseline values, generating subset profiles (\Cref{fig:subset}).
    
    \item We construct a slope vector field (\Cref{fig:vector}) where arrows indicate the direction of local $\beta_1$-$\beta_2$ correlations.
    
    \item For transverse curve parametrization, we introduce a new parameter $s$ and define an analytic curve $\tau$ transverse to the vector field (\Cref{eq:beta12}), marking baseline points $P=(\beta_1^*, \beta_2^*)$ from \Cref{tab:baseline}.
    
    \item We compute the profile likelihood for $c, D, s$, identify the peak $s^*$ yielding maximum likelihood, and obtain the corresponding point $Q=(\beta_1^{**}, \beta_2^{**})$ through \Cref{eq:beta12}. Linear interpolation is then used to draw the curves beginning from $Q$, which illustrate the relationship between $\beta_1$ and $\beta_2$  (\Cref{fig:r1}, \Cref{fig:r2}, \Cref{fig:r3}).
\end{itemize}

The key insights for our data are that for Regions 1 and 2, the unique peak at $s^*$ generates a contour curve through $Q$ representing all optimal $(\beta_1,\beta_2)$ pairs. In Region 3, the flat likelihood profile for $s$ indicates no unique optimum.

\end{enumerate}

\section{Application to experimental FRAP datasets}
\label{application}
\begin{figure}[htbp]
    \centering
    \begin{minipage}[t]{0.485\linewidth}
        \centering
        \includegraphics[width=\linewidth]{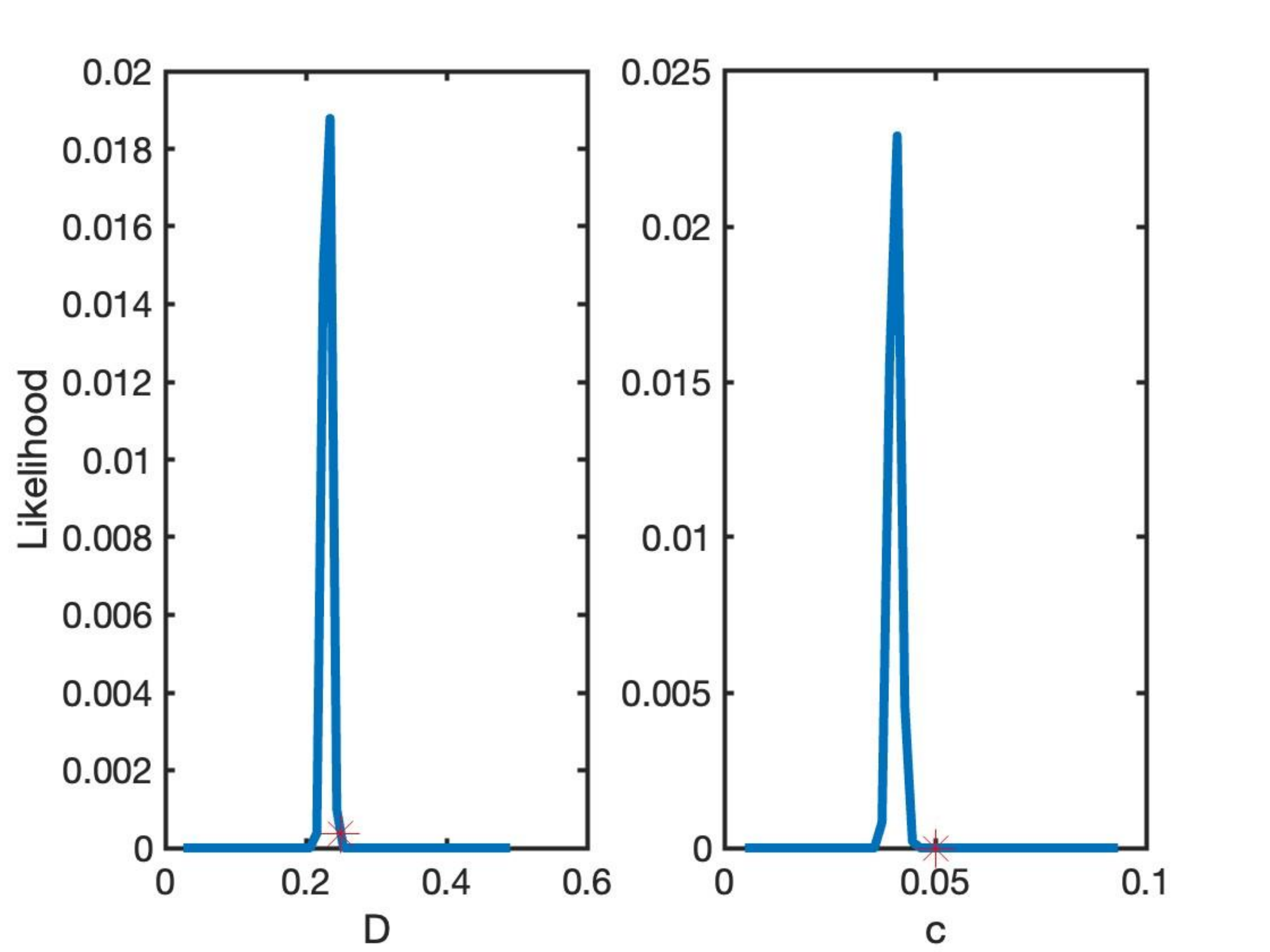}
    \end{minipage}
    \hspace{0.01\linewidth}
    \begin{minipage}[t]{0.485\linewidth}
        \centering
        \includegraphics[width=\linewidth]{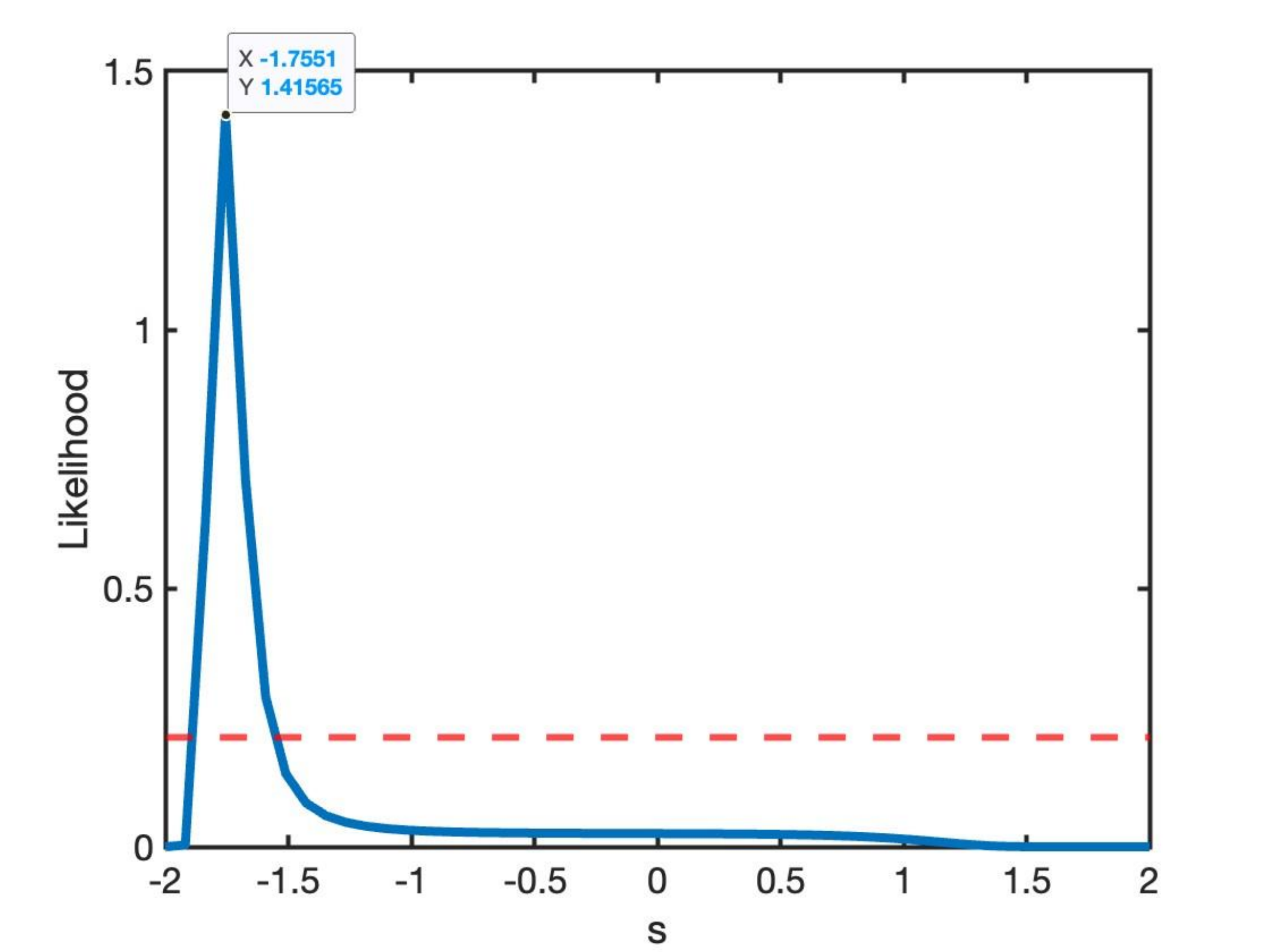}
    \end{minipage}
    
    \vspace{0.5em}
    \begin{minipage}[t]{0.7\linewidth}
        \centering
        \includegraphics[width=\linewidth]{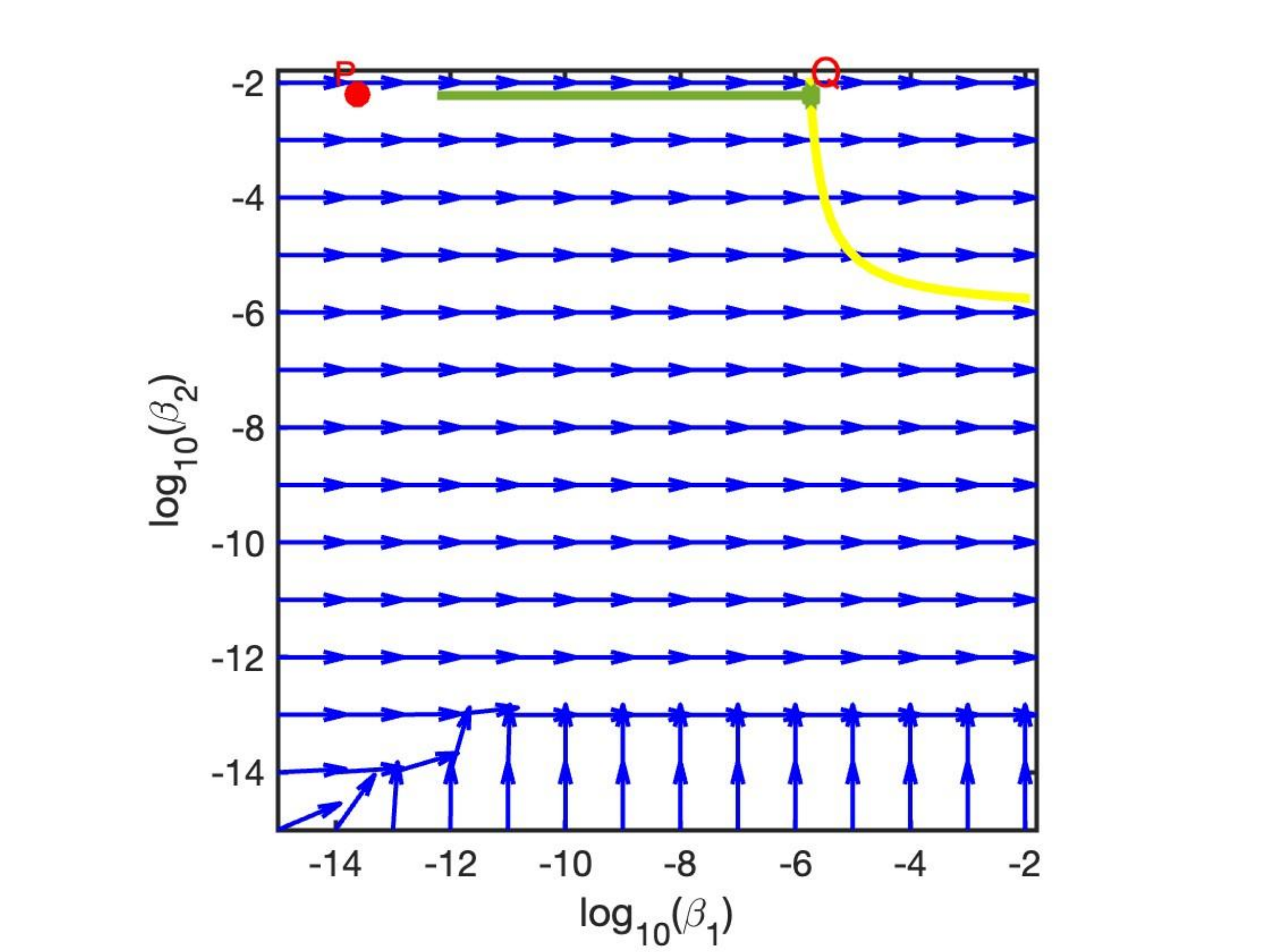}
    \end{minipage}

    \caption{Profile likelihoods for each interest parameter $c$, $D$ and $s$ in Region 1 given real experimental FRAP data in \Cref{fig:combined1} (top-left and top-right). $s$ achieves the maximum at $s^* = -1.7551$, corresponding to the point $Q$ on the yellow curve $\tau$, as well as the trace of error-minimizing green contour curve (bottom).}
    \label{fig:r1_real}
\end{figure}

\begin{figure}[htbp]
    \centering
    \begin{minipage}[t]{0.485\linewidth}
        \centering
        \includegraphics[width=\linewidth]{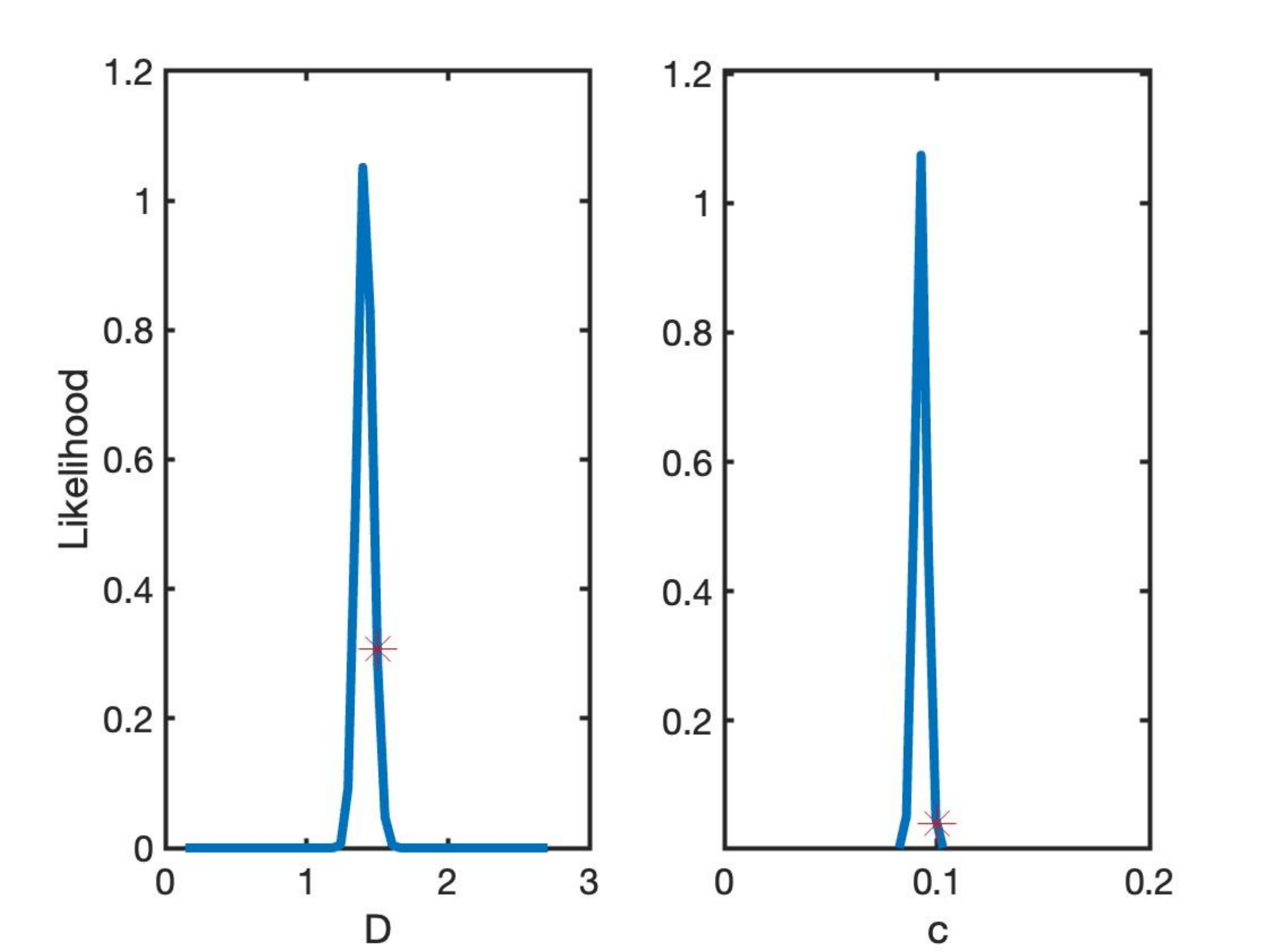}
    \end{minipage}
    \hspace{0.01\linewidth}
    \begin{minipage}[t]{0.485\linewidth}
        \centering
        \includegraphics[width=\linewidth]{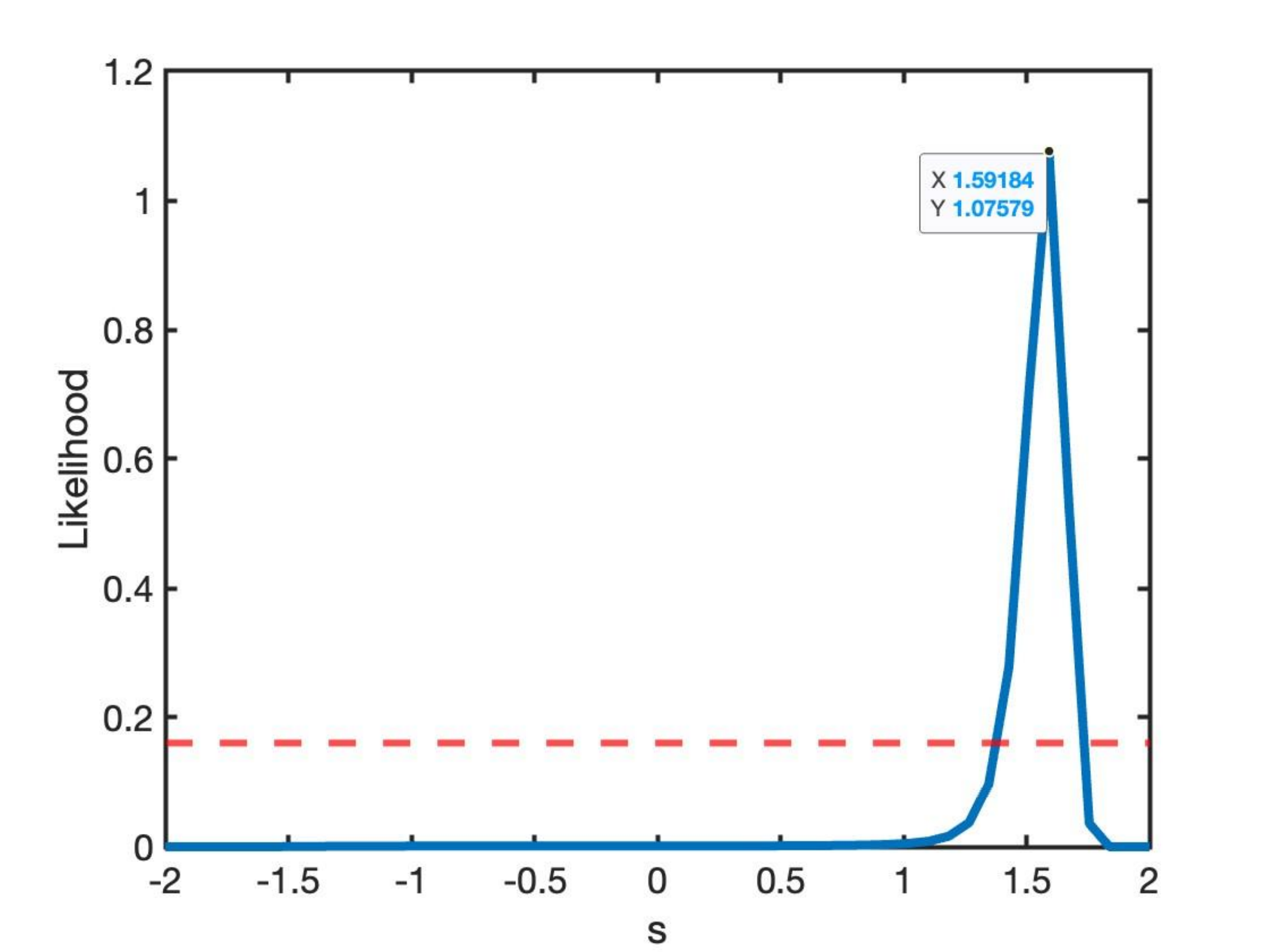}
    \end{minipage}
    
    \vspace{0.5em}
    \begin{minipage}[t]{0.7\linewidth}
        \centering
        \includegraphics[width=\linewidth]{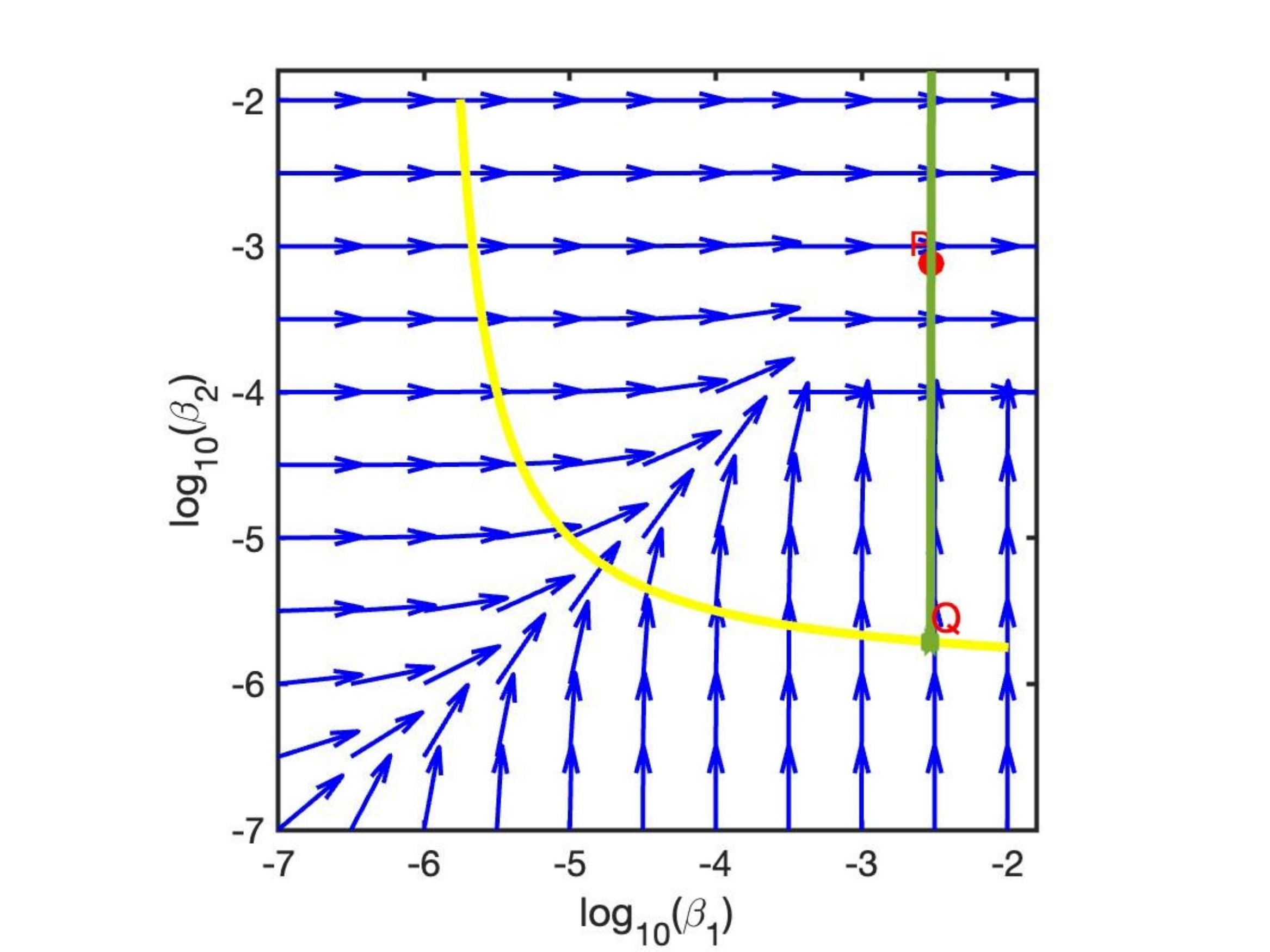}
    \end{minipage}

    \caption{Profile likelihoods for each interest parameter $c$, $D$ and $s$ in Region 2 given real experimental FRAP data in \Cref{fig:combined1} (top-left and top-right). $s$ achieves the maximum at $s^* = 1.59184$, corresponding to the point $Q$ on the yellow curve $\tau$, as well as the trace of error-minimizing green contour curve (bottom).}
    \label{fig:r2_real}
\end{figure}

\begin{figure}[htbp]
    \centering
    \begin{minipage}[t]{0.485\linewidth}
        \centering
        \includegraphics[width=\linewidth]{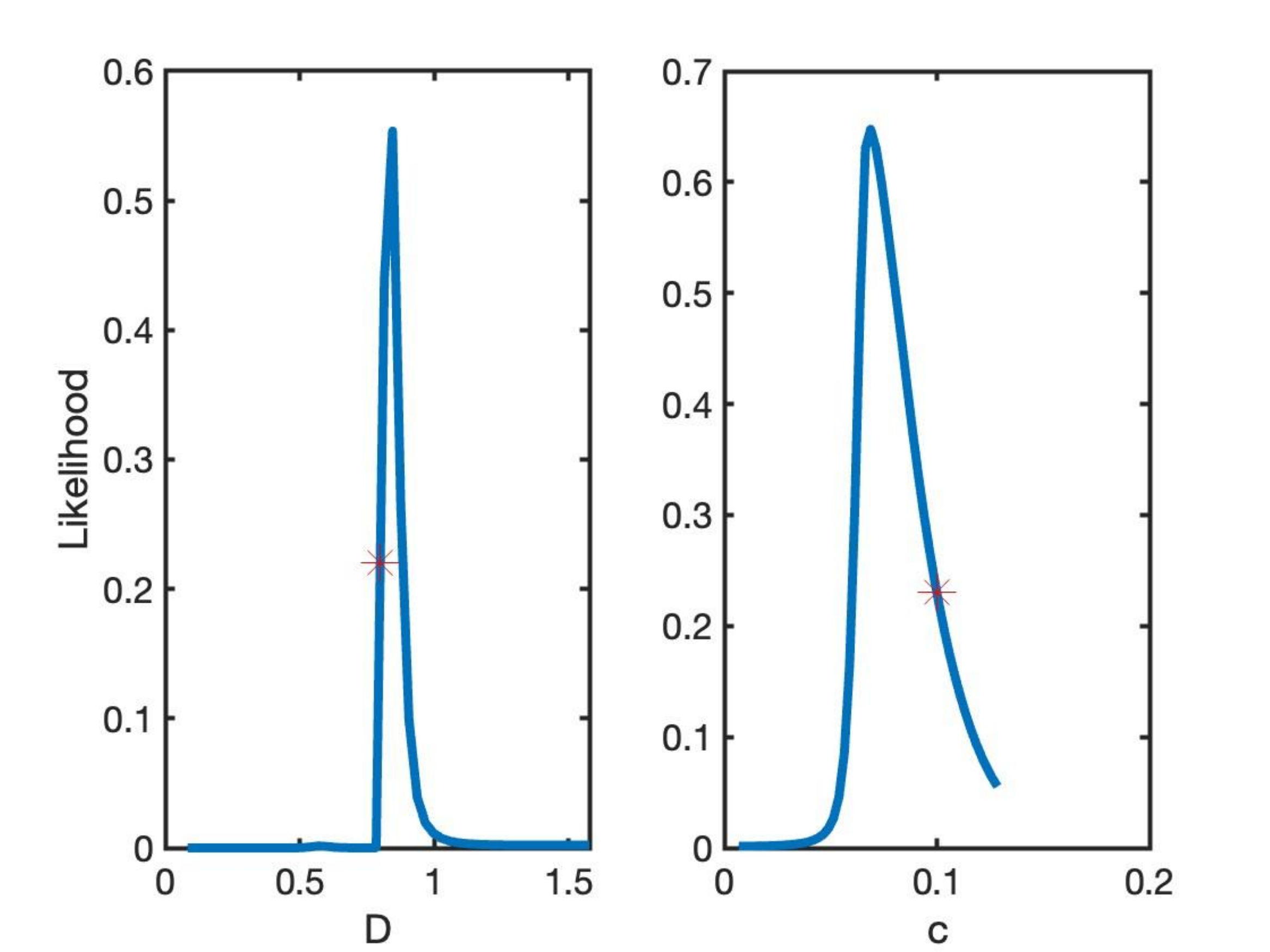}
    \end{minipage}
    \hspace{0.01\linewidth}
    \begin{minipage}[t]{0.485\linewidth}
        \centering
        \includegraphics[width=\linewidth]{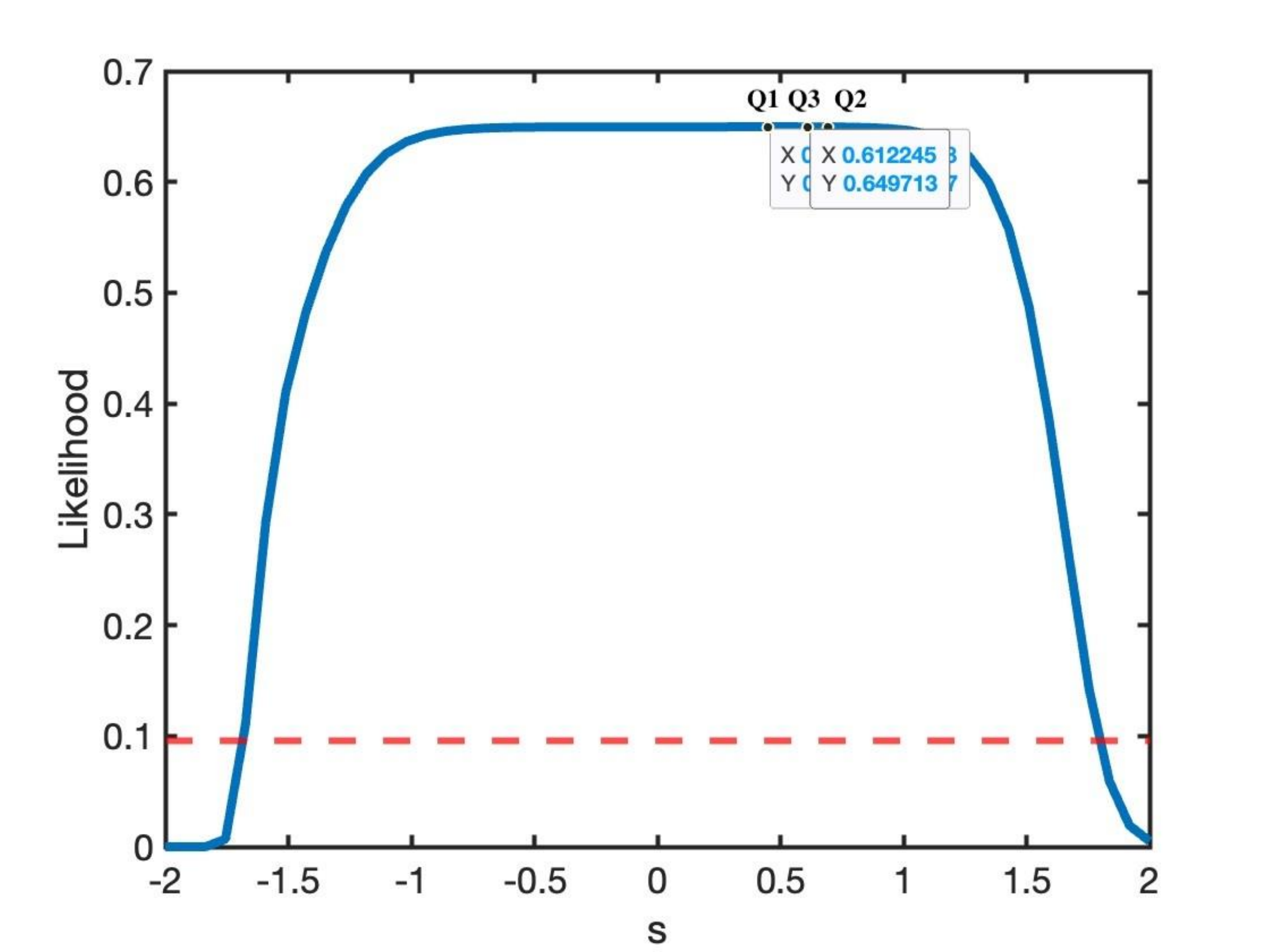}
    \end{minipage}
    
    \vspace{0.5em}
    \begin{minipage}[t]{0.7\linewidth}
        \centering
        \includegraphics[width=\linewidth]{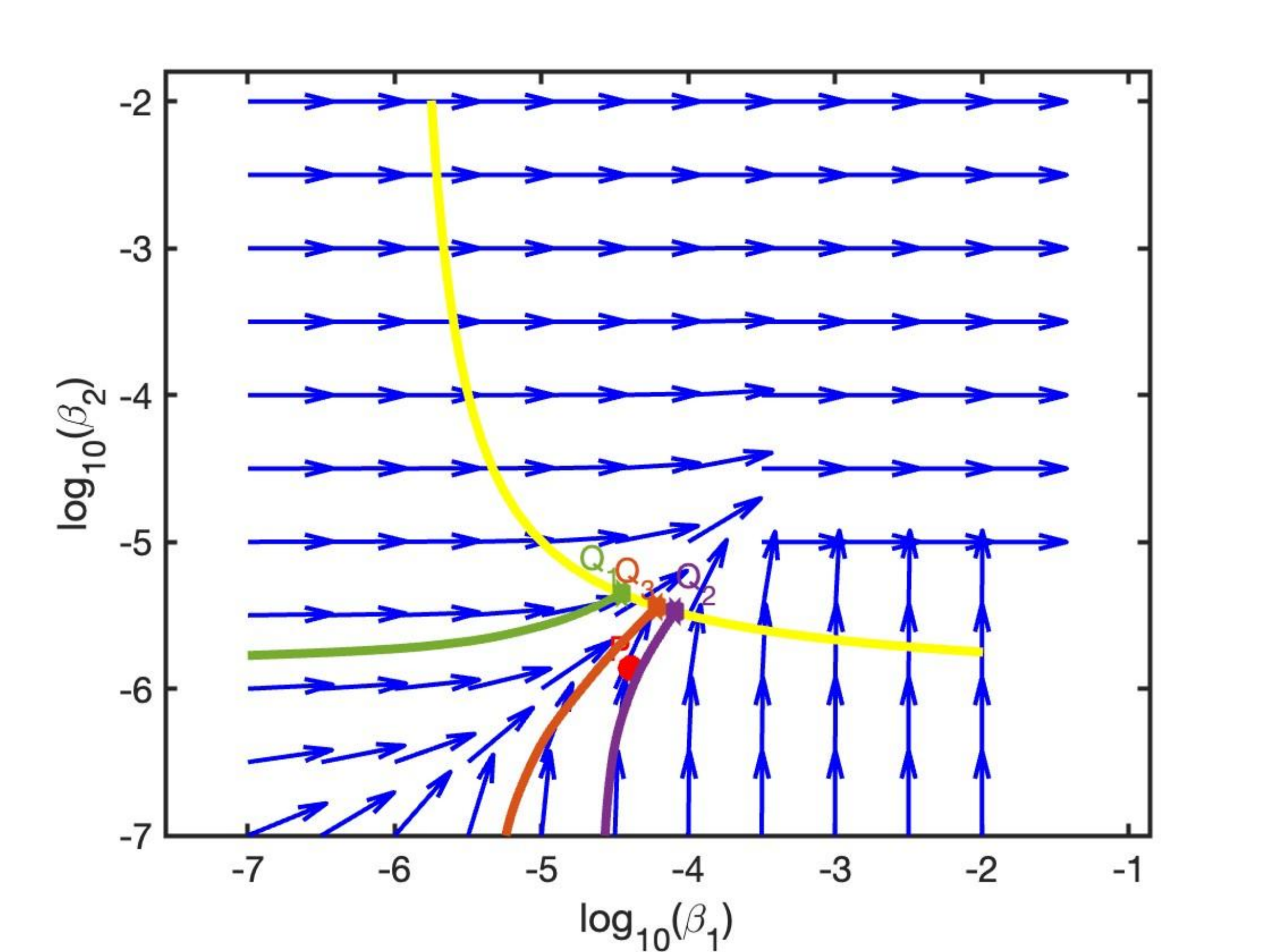}
    \end{minipage}

    \caption{Profile likelihoods for each interest parameter $c$, $D$ and $s$ in Region 3 given real experimental FRAP data in \Cref{fig:combined1} (top-left and top-right). Choosing three values for $s$: $s^* = 0.44898, 0.693878, 0.612245$ corresponds to points $Q_1, Q_2, Q_3$ on the yellow curve $\tau$, as well as the traces of error-minimizing contour curves (bottom).}
    \label{fig:r3_real}
\end{figure}

After validating this approach for synthetic FRAP datasets, we apply the outlined pipeline to the real average experimental FRAP data for each region. In particular, we generate the vector field and the corresponding contour curves using linear interpolation, but this time we use the real experimental FRAP data (blue points as shown in \Cref{fig:combined1}), and we set our baseline values as shown in \Cref{tab:real FRAP}. The pattern of the slope vector field and the contour curves in this case are given in \Cref{fig:r1_real}, \Cref{fig:r2_real} and \Cref{fig:r3_real}, which look very similar to those under our chosen baseline values, since the parameter values in \Cref{tab:real FRAP} are very similar to those in \Cref{tab:baseline}.

\section{Discussion}
Comparing the identifiability results that we obtained from profile likelihood (\Cref{fig:PL_all}) for $\beta_1$ and $\beta_2$ and the contour curves representing the relationship between $\beta_1$ and $\beta_2$ in slope vector fields (\Cref{fig:r1}, \Cref{fig:r2} and \Cref{fig:r3}), we observe the difference between practical non-identifiability and structural non-identifiability. For Region 1, it can be seen from the flat contour curve that although the value of $\beta_1$ is not unique, $\beta_2$ is fixed at $10^{-2}$. This is consistent with the likelihood profile graphs for $\beta_1$ and $\beta_2$, where $\beta_2$ is almost identifiable except that the right end of the curve does not go below the threshold. Similar conclusions also hold for Region 2, where $\beta_1$ is fixed at $10^{-3}$ for the main part of the green contour curve, coinciding with $\beta_1$ being more identifiable than $\beta_2$ as shown in the profile likelihood. On the other hand, for Region 3, we see from the profile likelihood curve that both $\beta_1$ and $\beta_2$ are structurally non-identifiable. Similarly, we find that the re-parameterized $s$ parameter is not identifiable in Region 3, making it difficult to determine a best contour curve representing the relationship between $\beta_1$ and $\beta_2$.

Our findings on parameter identifiability also carry important implications for the biological interpretation of RNA dynamics. In the context of FRAP experiments and PDE-based transport models, identifying which parameters can be uniquely estimated is crucial for linking quantitative model outputs to molecular mechanisms such as diffusion, active transport, and binding interactions. Being able to uniquely estimate parameters such as diffusion coefficient $D$ and transport velocity $c$ means that these values can be reliably connected to underlying RNA mobility and localization mechanisms. In contrast, the non-identifiability of binding and unbinding rates ($\beta_1, \beta_2$) indicates limitations of FRAP-based inference in estimating these parameters and suggests that additional experimental strategies may be needed. Overall, our analysis provides a rigorous foundation for using PDE models not only to fit experimental data but also to generate mechanistic insight into how RNAs are localized and regulated within cells.

\section*{Acknowledgments}
I would like to extend my heartfelt gratitude to my mentor, Professor Veronica Ciocanel, for her invaluable guidance, insightful feedback, and unwavering support throughout the course of my research.
\clearpage
\bibliographystyle{plain} 
\bibliography{my_references}

\end{document}